\documentclass[a4paper,12pt]{amsart2000}

\def\d{\operatorname{d}}
\def\ddiv{\operatorname{div}}
\def\lcf{\lbrack\! \lbrack}
\def\rcf{\rbrack\! \rbrack}
\def\nnabla{{\nabla\!\!\!\!\nabla}}

\def\Der{\operatorname{Der}}

\newtheorem{definition}{Definition}[section]
\newtheorem{lemma}[definition]{Lemma}
\newtheorem{theorem}[definition]{Theorem}
\newtheorem{proposition}[definition]{Proposition}
\newtheorem{corollary}[definition]{Corollary}

\begin{document}

\title{Divergence operators and odd Poisson brackets}

\author{Yvette Kosmann-Schwarzbach}
\address{
Centre de Math{\'e}matiques (U.M.R. 7640 du C.N.R.S.) \\
Ecole Polytechnique \\\mbox{F-91128} Palaiseau \\ France}
\email{yks@math.polytechnique.fr}

\author{Juan Monterde}
\address{Dpto. de Geometr{\'\i}a y Topolog{\'\i}a \\ Universitat de
Val{\`e}ncia \\
\mbox{E-46100} Burjasot (Val{\`e}ncia) \\ Spain}
\email{monterde@vm.ci.uv.es}

\date{}

\begin{abstract}
We define the divergence operators on a graded algebra, and we show
that, given an odd Poisson bracket on the algebra,
the operator that maps an element to the divergence of the hamiltonian
derivation that it defines is a generator of the bracket.
This is the ``odd laplacian'', $\Delta$,
of Batalin-Vilkovisky quantization.
We then study the generators of odd Poisson brackets on
supermanifolds,
where divergences of graded vector fields can
be defined either in terms of berezinian volumes or of
graded connections.
Examples include generators of the Schouten
bracket of multivectors on a manifold
(the supermanifold being the cotangent bundle where
the coordinates in the fibres are odd) and generators of the
Koszul-Schouten bracket of
forms on a Poisson manifold
(the supermanifold being the tangent bundle, with odd coordinates on
the fibres).
\end{abstract}


\maketitle

\section*{Introduction}

Graded algebras with an odd Poisson bracket~--~also called
Gerstenhaber
algebras~--~play an important role in the theory of deformations of
algebraic structures as well as in several areas of field theory, as
has been shown by Batalin and Vilkovisky \cite{B} \cite{BV}, Witten \cite{W},
Lian and Zuckerman
\cite{LZ}, Getzler \cite{G}, Hata and Zwiebach \cite{Z}, among others.
Generators of odd Poisson brackets, in the
sense of Equation \eqref{BV} below, 
are differential operators of order $2$
of the underlying graded algebra, sometimes called ``odd
laplacians'', and usually denoted by the letter~$\Delta$.
Batalin-Vilkovisky algebras~-- BV-algebras, 
for short~-- are a special class of these algebras, those for whose
bracket there exists a generator assumed to be of 
square~$0$. The geometrical approach to odd Poisson algebras and
BV-algebras in terms of supermanifolds 
was first developed by  Leites \cite{L}, Khudaverdian
\cite{kh0} (see also \cite{khn}) and Schwarz \cite{Sc}.

The purpose of this article is to study various constructions of
generators of odd Poisson brackets. 
Our constructions will rely on the general notion of divergence
operator on a graded algebra, which generalizes the
concept of the divergence of a vector field in elementary analysis.
Given an odd Poisson bracket, to each element in the algebra we
associate the divergence of the 
hamiltonian derivation that it defines, and 
we show that such a map from
the algebra to itself, multiplied by the factor ${\frac 12}$ and by an
appropriate sign,
is a generator of that bracket (Theorem \ref{generateur0}). 
We then adopt the language of supermanifolds to treat two
constructions which determine divergence operators on the
structural sheaf of the supermanifold.

The first construction uses berezinian volumes, and is modeled after
the construction of divergence operators on smooth, purely even 
manifolds which uses volume elements. 
Once a divergence operator is defined, we apply Theorem \ref{generateur0}
to obtain generators of an odd Poisson bracket on the supermanifold. 
One can deform any generator,
obtained
from a berezinian volume, by a change of berezinian volume, {\it i.e.},
the
multiplication by an invertible, even function. The deformed
generator then differs from the original one by
the addition of a hamiltonian derivation.
If the original generator is of square $0$, a sufficient
condition for a deformed generator to remain of square $0$
is given by a Maurer-Cartan equation. (Under the name ``BV quantum
master equation'', this condition plays a fundamental role in the
BV quantization of gauge fields.)
We study two special cases in detail: 
(i) the cotangent bundle of a
smooth manifold viewed as a supermanifold whose structural sheaf
is the sheaf of multivectors on the manifold, 
in which case the odd Poisson
bracket under consideration is the Schouten bracket, and
(ii) the tangent bundle of a smooth manifold as a
supermanifold whose structural sheaf is the sheaf of differential 
forms on the manifold, when the underlying smooth,
even manifold is a Poisson manifold.

The second construction utilizes graded linear connections on
supermanifolds, and generalizes 
to the graded case an approach to the
construction 
of divergence operators on purely even manifolds using linear
connections. Given a linear connection on a smooth, even
manifold, the divergence of a vector field $X$
is defined as the trace of 
the difference of the covariant derivation in the
direction of $X$ and of the map $[X,~.~]$, where $[~,~]$ is the Lie
bracket of vector fields.
In the case of a supermanifold, given a graded linear connection,
the divergence of a graded vector field is defined in a similar
manner, replacing the trace by the supertrace,
and the Lie bracket by the graded commutator.
We then apply Theorem \ref{generateur0}
to obtain a generator of an odd Poisson bracket on the supermanifold. 
As an example, we again study 
the generators of the Schouten bracket of multivectors on a manifold, 
equipped with a linear connection. 
On the one hand, there exists a unique 
generator of the Schouten bracket whose restriction to
the vector fields 
is the divergence defined by the linear connection.
On the other hand, on the cotangent bundle considered as a
supermanifold,
the linear connection on the manifold defines a
graded metric 
in a simple way. This graded metric in turn determines
a graded torsionless connection on this supermanifold -- the associated
Levi-Civita connection --, from which we obtain a generator of the Schouten
bracket, following our general procedure.
We show that these two generators of the Schouten bracket
coincide and that, 
when the connection is flat, this generator is of square $0$.

One can ask: what happens if we deal with an even Poisson bracket
instead of an odd one? The answer is that the phenomena 
in the odd and in the even cases are very different, although there is
a formal similarity of the constructions.
In the even case, the results extend those of the purely even case,
{\it e.g.}, the usual case of Poisson algebras of smooth
manifolds. 
Taking the divergence
of a hamiltonian vector field with respect to a volume form if the manifold
is orientable, or, more generally, to a density,  
yields a derivation of the algebra, 
{\it i.e.}, a vector field. 
This is the modular vector field that has been
studied in Poisson geometry (see \cite{Kz1}, \cite{We}) and in more
general contexts (see \cite{Hu2}, \cite{X}, \cite{KS}). 
One can prove that this vector field is
closed in the Poisson-Lichnerowicz cohomology.
A change in the volume element 
modifies the vector field by a hamiltonian vector field, therefore
the cohomology class of the closed vector field 
does not change. One thus obtains
a cohomology class, called ``the modular class'' \cite{We}.
So, while in the odd case we get a second-order differential operator
which is a generator of the bracket, in the even
case we get a first-order differential operator, in fact a derivation
of the structural sheaf of associative algebras.

\smallskip

The paper is organized as follows. In Section \ref{1}, we first 
recall the
definition of an odd Poisson bracket on a ${\mathbb Z}_2$-graded
algebra, and, in Section \ref{1.2},
we define algebraically the 
notion of a divergence
operator on a graded algebra and that of its curvature.
We then prove Theorem
\ref{generateur0}, which will serve as the main tool in our
constructions, and we study the deformation of divergence operators
and of the associated generators of odd Poisson brackets.

In Section \ref{2}, we study the divergence operators defined by
berezinian volumes (Proposition \ref{divberez}),
the associated generators of odd Poisson brackets, and their
deformation under a change of
berezinian volume. We show that the ``Batalin-Vilkovisky
quantum master equation''
appears as a sufficient condition for the modified
generator to remain of square 0
(Proposition \ref{master}).
In Section \ref{2.3},
we consider the example of the cotangent bundle of a manifold $M$,
considered as a supermanifold. To a volume element element $\mu$ on
$M$, there corresponds a berezinian volume, which behaves like the
``square of $\mu$''. In Theorem \ref{aa}, we show that
the generator of the Schouten bracket furnished by 
the general construction
outlined above
coincides with the generator obtained from the de~Rham
differential by the isomorphism defined by $\mu$
relating forms to multivectors.
In Section \ref{2.4}, we treat the case of the tangent bundle of $M$,
which is an odd Poisson supermanifold
when $M$ has a Poisson structure.
We prove that the generator defined by the divergence of hamiltonian
vector fields with respect to the canonical berezinian volume 
coincides with the Poisson homology operator, and that its square
therefore vanishes (Theorem \ref{bb}).
In Section \ref{2.5}, we express the properties of the supermanifolds
studied in Sections \ref{2.3} and \ref{2.4} 
in the language of QS, SP and QSP
manifolds of \cite{Sc} and \cite{ASZK} (Theorems \ref{piT*M} and \ref{piTM}).

In Section \ref{3}, we study the 
divergence operators defined by graded linear connections.
The definitions of a graded linear connection, its curvature and
torsion, and of the divergence operator that it defines (Proposition
\ref{propdivconn})
are given in Section \ref{divop}.
We then study the generator associated to a
torsionless graded linear connection, of an odd Poisson bracket 
and the effect of a change of
connection on the generator (Section \ref{3.2}).  
The Levi-Civita connection of
a graded metric on a supermanifold is introduced in Section \ref{3.3}.
In the remaining part of Section \ref{3},  we study the cotangent
bundle of a manifold as an odd Poisson supermanifold. More
specifically we study two constructions of generators 
of the Schouten bracket
associated to a torsionless linear connection on the base manifolds,
and we show that the two constructions yield the same generator 
(Theorem \ref{cc}).
We conclude the paper with remarks concerning the relationships
between divergence operators, right and left 
module structures in the theory of
Lie-Rinehart algebras \cite{Hu} and right 
and left ${\mathcal D}$-modules, and their analogues in the graded
case \cite{P}, and we 
formulate a conjecture 
regarding the existence of a unique prolongation of a
divergence operator on a supermanifold into a generator of an odd
bracket on the algebra of graded multivectors.

\smallskip

We shall usually denote a supermanifold by a pair
$(M,\mathcal A)$ where $M$ is an ordinary smooth manifold,
called the base
manifold, and $\mathcal A$ is a sheaf over $M$ of ${{\mathbb Z}}_{2}$-graded
commutative, associative algebras.
The sections of $\mathcal A$ will be
denoted by $f,g, \dots$, but this notation will be modified in some
instances. When $a$ is an element of a ${{\mathbb Z}}_{2}$-graded vector space,
$|a|$ denotes the ${{\mathbb Z}}_{2}$-degree of $a$
and, whenever it appears in
a formula, it is understood that $a$ is homogeneous. The word
``graded'' will often be omitted.
The bracket $[~,~]$ denotes the graded commutator.
Manifolds and maps are assumed to be smooth.
We recall some general properties of
supermanifolds and the definition of the berezinian
volumes in the Appendix.

\section{Odd Poisson
brackets and divergence operators} \label{1}

In this section, we review the main definitions concerning odd
Poisson brackets on graded algebras and we study how to construct
generators of such brackets.

\subsection{Gerstenhaber and BV-algebras}
Let ${\bf A}$ be a ${{\mathbb Z}}_2$-graded commutative, associative
algebra over a field of characteristic $0$. 
The multiplication 
map of algebra ${\bf A}$, \ $(f,g) \in {\bf A} \times {\bf A} 
\mapsto fg \in {\bf A}$, will be
denoted by $m$. By definition, an {\it odd Poisson bracket} or a
${{\mathbb Z}}_2$-{\it Gerstenhaber bracket} on ${\bf A}$
is an odd bilinear map, $\pi~: (f,g) \in {\bf A} \times {\bf A} \mapsto
\lcf f ,g \rcf \in {\bf  A}$, satisfying,
for any $f,g,h \in {\bf A}$,
\begin{itemize}
\item $\lcf f,g\rcf = -(-1)^{(|f|-1)(|g|-1)}\lcf g,f\rcf$ \quad
(skew-symmetry),
\item $\lcf f,\lcf g,h\rcf\rcf =
 \lcf\lcf f,g\rcf,h\rcf  + (-1)^{(|f|-1)(|g|-1)}
 \lcf g,\lcf f,h\rcf\rcf$ \quad (graded Jacobi identity),
\item $\lcf f,gh\rcf = \lcf f,g\rcf h +(-1)^{(|f|-1)|g|}g \lcf f,h\rcf$
\quad (Leibniz rule) \ .
\end{itemize}
(``Map $\pi$ is odd'' means that $|\lcf f,g\rcf| = |f|+|g|-1$  modulo $2$.)
The pair $({\bf A}, \pi)$ is then called an
{\it odd Poisson
algebra} or a ${{\mathbb Z}}_2$-{\it Gerstenhaber algebra}.

A linear map  of odd degree, $\Delta : {\bf A} \rightarrow {\bf A}$,
such that, for all $f, g \in {\bf A}$,
\begin{equation} \label{BV}
\lcf f,g  \rcf = (-1)^{|f|} \left(\Delta(fg) -
(\Delta f) g - (-1)^{|f|}f(\Delta g)\right) \ ,
\end{equation}
is called a {\it generator} or a {\it generating operator}
of $\pi$ (or of bracket $\lcf ~,~ \rcf$).
If there exists a generator $\Delta$ of the bracket which is
of square $0$, then $({\bf A}, \pi, \Delta)$ is 
called a
${{\mathbb Z}}_2$-{\it Batalin-Vilkovisky algebra}, or 
{\it BV-algebra} for short.

\medskip

\noindent{\it Remark.} 
Since a Gerstenhaber bracket in the usual sense, 
defined on a ${{\mathbb Z}}$-graded algebra,
is of ${{\mathbb Z}}$-degree $-1$, it is clear that it
can also be considered to be
a ${{\mathbb Z}}_2$-Gerstenhaber
bracket. Similarly, Batalin-Vilkovisky algebras in the usual, 
${{\mathbb Z}}$-graded sense, are particular cases of  
${{\mathbb Z}}_2$-Batalin-Vilkovisky algebras.

\medskip

A generator of an odd Poisson bracket is clearly not a derivation of
the graded associative algebra $({\bf A},m)$, unless the bracket is
identically $0$. We shall now see under what condition a generator
$\Delta$ of an
odd Poisson bracket $\pi$ is a derivation of the graded Lie algebra
$({\bf A},\pi)$. A straightforward computation,
using the
defining relation \eqref{BV} for a generator, yields the identity,
\begin{equation} \label{ww}
\Delta^{2}(fg) -(\Delta^{2}f)g - f(\Delta^{2}g) =
(-1)^{|f|}\left(\Delta(\lcf
f,g \rcf)
- \lcf \Delta f ,g \rcf - (-1)^{|f|-1}\lcf f, \Delta g \rcf \right) \ ,
\end{equation}
for all $f,g \in {\bf A}$. 
From \eqref{ww}, we obtain
\begin{lemma} \label{ccc}
A generator 
$\Delta$ of an odd Poisson bracket
$\pi$
is an odd derivation of the odd
Poisson algebra, $({\bf A}, \pi)$, if and only if 
the map $\Delta^{2}$
is an even derivation of the graded associative
algebra
$({\bf A}, m)$. In particular, if $\Delta^2 = 0$, then $\Delta$ is a
derivation of $({\bf A}, \pi)$.
\end{lemma}

\noindent{\it Remark.} In the case of a ${\mathbb Z}$-graded
algebra and a the linear map $\Delta$ of degree
$-1$, the map $\Delta ^{2}$, which is of degree $-2$, is a
derivation of $({\bf A},m)$ if and only if it vanishes.
Therefore, in this case, a generator $\Delta$ of $\pi$
is a derivation of $({\bf A},\pi)$ if and only if $\Delta ^2=0$.

\medskip

The following identity, of which we
shall make use in Section  \ref{2.5}, is the result of another computation. 
If $D$ is an odd derivation 
of $({\bf A},m)$, then
\begin{equation} \label{derdelta}
[D,\Delta](fg)- ([D,\Delta]f)g - f([D,\Delta]g) = 
(-1)^{|f|} 
\left( D\lcf f,g \rcf - \lcf Df,g \rcf - (-1)^{|f|-1}\lcf f,Dg \rcf \right).
\end{equation}

Let $\delta_m$ be the graded Hochschild
differential of algebra $({\bf A},m)$. Equation \eqref{BV} expresses 
the equality
$$\pi (f,g) = (-1)^{|f|-1} (\delta_m \Delta)(f,g) \ .$$
Let us also introduce the graded Chevalley-Eilenberg differential of
$({\bf A}, \pi)$, denoted by $\delta_{\pi}$. Equation \eqref{ww}
expresses the equality 
$$(\delta_{\pi} \Delta)(f,g) = (-1)^{|f|-1} (\delta_m (\Delta^2))(f,g) \ ,$$
and Lemma \ref{ccc} can be reformulated as follows~: $\Delta$ is a 
$\delta_{\pi}$-cocycle 
if and only if $\Delta^2$ is a 
$\delta_{m}$-cocycle.

\subsection{Divergence operators} \label{1.2}

Let $({\bf A},m)$ be a ${{\mathbb Z}}_2$-graded commutative, associative
algebra, and let $\Der {\bf A}$ be the graded vector space of graded
derivations of $({\bf A},m)$. 
By definition, a {\it divergence operator} on ${\bf A}$ 
is an even linear map, $\ddiv : \Der{\bf A} \to  {\bf A}$, such that
\begin{equation} \label{gooddiv} 
\ddiv(fD) = f\ \ddiv(D) + (-1)^{|f||D|}D(f) \ ,
\end{equation} 
for any $D \in \Der{\bf A}$ and any $f \in {\bf A}$.

This definition obviously generalizes the usual notion of divergence
of vector fields in elementary analysis. More generally, if ${\bf A}$
is the purely even algebra of smooth real- or complex-valued functions
on a smooth manifold, divergence operators on ${\bf A}$ 
can be defined by means of
either volume forms or linear connections, the two approaches being
related in a simple way. See \cite{Kz1}. These two approaches
will be generalized below in Sections \ref{2} and \ref{3}, respectively.

We define the {\it curvature
of a divergence operator},
as the bilinear map,
${\mathcal R}^{\ddiv}~: \Der{\bf A} \times \Der{\bf A} \to {\bf A}$
by
\begin{equation} \label{curvdiv} 
{\mathcal R}^{\ddiv} (D_1,D_2)
=\ddiv [D_1,D_2]-D_1(\ddiv D_2) + (-1)^{|D_1||D_2|}D_2(\ddiv D_1) \ .
\end{equation} 
for any derivations $D_1, D_2 \in \Der {\bf A}$.
A short computation shows that  ${\mathcal R}^{\ddiv}$ 
is ${\bf A}$-bilinear.
If $\delta_{[~,~]}$ denotes the graded Chevalley-Eilenberg differential
of the
graded Lie algebra $(\Der{\bf A}, [~,~])$ 
acting on cochains on $\Der{\bf A}$
with values in the $\Der{\bf A}$-module ${\bf A}$, then the
definition of $\mathcal R^{\ddiv}$ in \eqref{curvdiv} can be written
$$
{\mathcal R}^{\ddiv} = \delta_{[~,~]}(\ddiv) \ .
$$ 
  
\subsection{Divergence operators and generators} \label{1.3}
On the odd Poisson algebra $({\bf A}, \pi)$, we
define the hamiltonian mapping by
$X^{\pi}: {\bf A} \to \Der{\bf A}$, defined by
$$
f \mapsto X^{\pi}_f = \lcf f,~.~\rcf \ .
$$
The graded Jacobi identity for the bracket $\pi$ is equivalent to 
the relation
\begin{equation} \label{hamilt}
X^{\pi}_{\lcf f, g \rcf} = [X^{\pi}_f,X^{\pi}_g] \ ,
\end{equation} 
for any $f, g \in {\bf A}$, which expresses the fact that $X^{\pi}$ is a
morphism of graded Lie algebras from $({\bf A}, \pi)$ to 
$(\Der{\bf A}, [~,~])$.

We now introduce the main object of interest in this paper, the
odd linear map $\Delta~:{\bf A} \to {\bf A} $, 
depending on both $\pi$ and the choice of a divergence operator,
defined by
\begin{equation} \label{generator}
\Delta f= (-1)^{|f|}\frac 12 \ddiv(X^{\pi}_f) \ ,
\end{equation} 
for $f \in {\bf A}$.

\begin{theorem} \label{generateur0}
The operator $\Delta$ on ${\bf A}$, defined by \eqref{generator},
is a generator of bracket
$\pi$.
\end{theorem}
\noindent{\it Proof.} To show that Equation \eqref{BV} is satisfied,
we compute $\Delta (fg)$ using the Leibniz rule for the odd
Poisson bracket, and the fundamental property \eqref{gooddiv} of the
divergence operators,
$$\aligned
(-1)&^{|f|+|g|}\Delta (fg)
= \frac 12 \ddiv X^{\pi}_{fg}
= \frac 12 \ddiv (f X^{\pi}_{g}+
(-1)^{|f||g|}g X^{\pi}_{f}) \\
 &=
(-1)^{|g|}f \Delta g +
(-1)^{|f|(|g|+1)}g \Delta f
 +
\frac 12 (-1)^ {|f|(|g|+1)} \lcf g,f\rcf +
\frac 12 (-1)^ {|g|} \lcf f,g\rcf
\\
 &=
(-1)^{|g|}f \Delta g +
(-1)^{|f|+|g|}(\Delta f)g +
 (-1)^ {|g|} \lcf f,g\rcf \ ,
\endaligned
$$
and the result follows.\qed

This fundamental result must be contrasted with a parallel, but
strikingly different result valid for even Poisson brackets, and
consequently also, in the usual case, for ungraded Poisson algebras.
If, in \eqref{generator}, we replace the hamiltonian operator defined
by an odd Poisson bracket by the one defined by
an even Poisson bracket, a computation similar to
the proof of Theorem \ref{generateur0} 
shows that the operator thus defined is a
derivation of the associative multiplication. When the Poisson
algebra is the
algebra of functions of a smooth orientable
Poisson manifold, in which case the divergence
operator is the one associated with a volume element, this
derivation is a vector field. It is easy to see that, 
up to the factor $1/2$, it
coincides with the {\it modular vector field}  of the Poisson manifold
\cite{We} (also \cite{Hu2}, \cite{KS} and references cited therein).

\medskip

We now establish a relation between the operator $\Delta$ defined by
$\pi$ and $\ddiv$, and the curvature of the
divergence operator evaluated on hamiltonian derivations.

\begin{proposition} For any $f, g \in {\bf A}$,
\begin{equation} \label{eqcurv}
\Delta \lcf f,g \rcf - 
\lcf \Delta f ,g \rcf - (-1)^{|f|-1}\lcf f, \Delta g \rcf   
= (-1)^{|f|+|g|-1} \frac {1}{2} 
{\mathcal R}^{\ddiv}(X^{\pi}_f,X^{\pi}_g) \ .
\end{equation}
\end{proposition}
\noindent{\it Proof.}
By \eqref{hamilt} and the definition of $\Delta$,
$$
\begin{aligned}
\Delta \lcf f,g\rcf &= 
(-1)^{|f|+|g|-1}{\frac 12} \ddiv X^{\pi}_{\lcf f,g\rcf} 
= (-1)^{|f|+|g|-1} {\frac 12} \ddiv
[X^{\pi}_{f},X^{\pi}_{g}] \ , \\
\lcf f,\Delta g \rcf &=
(-1)^{|g|}\frac 12 X^{\pi}_{f}(\ddiv X^{\pi}_{g}) \ , \qquad
\lcf \Delta f, g \rcf =
(-1)^{|f||g|-1}\frac 12 X^{\pi}_{g}(\ddiv X^{\pi}_{f}) \ .
\end{aligned}
$$
In view of the definition of the
curvature, the proposition follows.
\qed

\begin{corollary} \label{dercurv}
The generator $\Delta$ is a derivation of $({\bf A},\pi)$
if and only if ${\mathcal R}^{\ddiv}$ vanishes on the hamiltonian
derivations.
\end{corollary}

In terms of the differentials $\delta_{\pi}$ 
and  $\delta_{[~,~]}$, Equation \eqref{eqcurv} can
be written
$$
(\delta_{\pi}\Delta)(f,g) = (-1)^{|f|+|g|-1} {\frac {1}{2}}(\delta_{[~,~]} \ddiv)
(X^{\pi}_f,X^{\pi}_g) \ .
$$

\subsection{Deformations of divergence operators and of generators}
Since the difference of 
two divergence operators is an ${\bf A}$-linear map
from $\Der{\bf A}$ to ${\bf A}$, the space of divergence operators
on ${\bf A}$ is an affine space over 
${\rm Hom}_{\bf A}(\Der{\bf A},{\bf A})$.
We shall be interested in the case where 
the difference of two divergence operators is
an evaluation map, $D \in \Der{\bf A} \mapsto D(2w) \in {\bf A}$,
where $w$ is a fixed, even element in $\bf A$. (The factor $2$ is
conventional.)

Since the difference of two generators of an odd Poisson bracket
$\pi$ on ${\bf A}$ is a derivation of $({\bf A},m)$, the space of
generators
of $\pi$ is an affine space over $\Der{\bf A}$.

\begin{proposition}\label{prime}
Let $\ddiv$ and $\ddiv ^{\prime}$
be divergence operators on $\bf A$ 
such that there exists an even $w \in {\bf A}$ satisfying  
\begin{equation} \label{deformdiv}
\ddiv ^{\prime} D = \ddiv D +D(2w) \ ,
\end{equation} 
for all
$D \in \Der{\bf A}$.
For $\pi$ a fixed odd Poisson bracket, let $\Delta$ and 
$\Delta ^{\prime}$ be the
generators of the bracket defined by \eqref{generator} for $\ddiv$ and 
$\ddiv^{\prime}$, respectively. Then
\begin{equation} \label{defgen}
\Delta ^{\prime} = \Delta + X^{\pi}_w \ .
\end{equation} 
\end{proposition}
\noindent{\it Proof.}  Equation \eqref{defgen} follows from \eqref{deformdiv} and the skew-symmetry of the bracket.
\qed

\medskip

\noindent{\it Remark.}
In the case of an even Poisson algebra, and, in particular, in the
usual case of the Poisson algebra of a smooth manifold, a similar
argument is valid.
As a consequence, one proves that the class of the modular
vector field in the Poisson cohomology is well-defined, independently
of the choice of a volume element \cite{We}.

\medskip

We shall now consider under what condition a generator
of $\pi$  with vanishing square 
remains of square $0$ when modified by the
addition of an interior derivation $X^{\pi}_w$.

\begin{proposition} \label{genprime}
If $\Delta$ is a generator of square $0$ of bracket $\pi$, and $w$ is
an even element of ${\bf A}$, then
\begin{equation} \label{777}
(\Delta + X^{\pi}_w)^2 = X^{\pi}_{\Delta w + \frac {1}{2}
\lcf w, w \rcf} \ .
\end{equation}
\end{proposition}
\noindent{\it Proof.} From Lemma \ref{ccc}, we know that 
$\Delta ^2 =0$ implies that $\Delta$ is a derivation of $\pi$, whence
$[\Delta, X^{\pi}_w] = X^{\pi}_{\Delta w}$.
Using this relation and \eqref{hamilt}, we obtain
$$
(\Delta + X^{\pi}_w)^2
= [\Delta, X^{\pi}_w] + \frac {1}{2} [X^{\pi}_w,X^{\pi}_w]
= X^{\pi}_{\Delta w}+ \frac{1}{2} X^{\pi}_{\lcf w,w \rcf} \ ,
$$
whence the result.
\qed

The equation
\begin{equation} \label{MC}
\Delta w + \frac {1}{2} \lcf w,w \rcf= 0 \ ,
\end{equation} 
is the
{\it Maurer-Cartan equation}, familiar from deformation
theory. (See \cite{St}.) Thus we can state 

\begin{corollary}
Let $\Delta$ be a generator of square $0$ of bracket $\pi$, and let
$w$ be an even element of ${\bf A}$.
The generator 
$\Delta ^{\prime}= \Delta + X^{\pi}_w$ is of square $0$ if and
only if the hamiltonian operator $X^{\pi}_{\Delta w + \frac {1}{2}
\lcf w, w \rcf}$
vanishes.
If, in particular, 
$w$ satisfies the Maurer-Cartan equation \eqref{MC}, then 
the generator 
$\Delta ^{\prime}= \Delta + X^{\pi}_w$ is of square $0$.
\end{corollary}

\section{Berezinian volumes and generators of odd Poisson
brackets} \label{2}

An {\it odd Poisson supermanifold} (resp., a {\it BV-supermanifold}) 
is a supermanifold
$(M,{\mathcal A})$ whose sheaf of
functions, ${\mathcal A}$, is a sheaf of odd Poisson algebras (resp., of
BV-algebras). 
In the context of supermanifold theory, an odd Poisson bracket
is often referred to as
an {\it antibracket} \cite{B} \cite{W} \cite{Z} or a
{\it Buttin bracket} \cite{L}. The notions of derivations and
divergence operators that we have introduced in
Section \ref{1} have
obvious analogues in the case of sheaves of algebras, and we shall use
the same symbols.  
On a supermanifold, a derivation of the sheaf of functions is called a
{\it graded vector field}.
We use the term ``operator from ${\mathcal A_1}$ to ${\mathcal A}_2$'' 
for a
morphism of sheaves of vector spaces
from ${\mathcal A}_1$ to ${\mathcal A}_2$.

In this
section, we show how generators of an odd
Poisson bracket on a supermanifold can be obtained from berezinian 
volumes.

\subsection{Divergence operators defined by berezinian volumes}

We first recall the main properties of the Lie derivatives of berezinian
sections. (See the Appendix for the definition of the berezinian
sheaf. See,
{\it e.g.}, \cite{D} for a proof of the following proposition.)

\begin{proposition} \label{lieder}
Let $\xi$ be a berezinian section on $(M,{\mathcal A})$. For any
graded vector field  $D$ and any
section $f$ of ${\mathcal A}$, 
\begin{equation} \label{lieder1}
{\mathcal L}_D(\xi . f)=
{\mathcal L}_D(\xi).f + (-1)^{|D||\xi|} \xi.D(f) \ ,
\end{equation}
and
\begin{equation}\label{lieder2}
{\mathcal L}_{f.D}(\xi) = (-1)^{|f|(|D|+|\xi|)}{\mathcal L}_D(\xi . f) \ .
\end{equation}
\end{proposition}

We shall now define the divergence operator associated with a
berezinian volume.. 
\begin{proposition} \label{divberez}
Let $\xi$ be a berezinian volume. 
For any graded vector field $D$,
there exists a unique section, $\ddiv_\xi(D)$, of $\mathcal A$ such that
\begin{equation} \label{divxi}
{\mathcal L}_D\xi = (-1)^{|D||\xi|} \xi.\ddiv_\xi(D) \ .
\end{equation}
The map $D \mapsto
\ddiv_\xi(D)$ from ${\rm \Der}{\mathcal A}$ to ${\mathcal A}$ 
defined by \eqref{divxi} is a divergence operator.
\end{proposition}

\noindent{\it Proof.}
The map $\ddiv_\xi$ 
is even, since
$|\xi|+|\ddiv_\xi(D)| =|\xi.\ddiv_\xi(D)| = |{\mathcal L}_D\xi| =
|D|+|\xi|.$
We must prove that 
$\ddiv_\xi$ satisfies \eqref{gooddiv}.
Using Proposition \ref{lieder}, we obtain
$$\aligned
\xi.\ddiv_\xi(fD) &=  (-1)^{(|f|+|D|)|\xi|} {\mathcal L}_{fD}\xi
= (-1)^{|D|(|f|+|\xi|)}{\mathcal L}_D(\xi.f)\\
 &= (-1)^{|D|(|f|+|\xi|)}({\mathcal L}_D(\xi).f + (-1)^{|D||\xi|} 
\xi.D(f)) \\
&= (-1)^{|f||D|}(\xi.\ddiv_\xi(D)f + \xi.D(f))  \ ,
\endaligned
$$
whence the result. \qed

\medskip

\noindent{\sc Example.}
If $(M,{\mathcal A}) = {\mathbb R}^{m|n}$ with graded coordinates
$(x^{1},\dots,x^{m},s^{1},\dots,s^{n})$, the section
$$\xi =
[\d^Gx^{1}\wedge\dots\wedge\d^Gx^{m} \otimes
\frac\partial{\partial s^{1}}\circ\dots\circ\frac\partial
{\partial s^{n}}]
$$
is a berezinian volume and, if
$D = \sum_{i=1}^m g^{i} {\frac  {\partial}{\partial x^{i}}} +
\sum_{\rho = 1}^n h^{\rho}\frac {\partial}{\partial s^{\rho}}$, then
$$\ddiv_\xi(D) = \sum_{i=1}^m \frac {\partial g^{i}}{\partial x^{i}} +
\sum_{\rho= 1}^n (-1)^{|h^{\rho}|}
\frac {\partial h^{\rho}}{\partial s^{\rho}} \ .
$$

\medskip

A short computation shows that, for graded vector fields $D_1$ and
$D_2$,
\begin{equation}\label{commutator}
\ddiv_{\xi} [D_1,D_2] = D_1(\ddiv_{\xi}D_2) 
- (-1)^{|D_1| |D_2|} D_2(\ddiv_{\xi}D_1) \ .
\end{equation}
In view of the definition of the curvature of a divergence operator 
\eqref{curvdiv}, we have proved
\begin{proposition} \label{curvber}
For any berezinian volume $\xi$, 
the curvature ${\mathcal R}^{\ddiv_{\xi}}$ 
of the divergence operator $\ddiv_{\xi}$ vanishes.
\end {proposition}

We now consider the effect on the divergence operator
of a change of berezinian volume.
When $v$ is an invertible, even section of $\mathcal A$, then
$\xi.v$ is also a generator of the berezinian sheaf.
We remark that any invertible, even section,
$v$ of ${\mathcal A}$, can be written as
$\pm e^{2w}$ for an even section $w$ of ${\mathcal A}$.
In fact, $v$ can be written as the product of a nowhere vanishing
function on the base manifold and a function $1+u$, where $u$ is
nilpotent. Since $u$ is nilpotent, say of order $k$,
$1+u$ is equal to ${\rm exp}({\rm ln}(1+u))$, where
${\rm ln}(1+u)= \sum_{j=1}^{k} (-1)^{j-1} {\frac {u^{j}}{j}}$.

\begin{proposition} \label{newvolume}
Let $\xi$ be a berezinian volume.
For any invertible, even section $v$ of $A$, the berezinian section
$\xi.v$ is a berezinian volume,
and, for any graded vector field $D$,
\begin{equation} \label{conformalchange}
\ddiv_{\xi.v}(D) = \ddiv_\xi(D) + v^{-1}D(v) \ .
\end{equation}
If $v = e^{2w}$, where $w$ is an even section of ${\mathcal
A}$,
then
\begin{equation} \label{particular}
\ddiv_{\xi.e^{2w}}(D) = \ddiv_\xi(D) + D(2w) \ .
\end{equation}
\end{proposition}

\noindent{\it Proof.}
Using \eqref{lieder1}, we compute
$$
\aligned
(\xi.v).\ddiv_{\xi.v}(D)&= (-1)^{|D||\xi|} {\mathcal L}_D(\xi.v)\\
&= (-1)^{|D||\xi|} {\mathcal L}_D(\xi).v+\xi.D(v)\\
&=  \xi.(\ddiv_\xi(D)v)+\xi.D(v) \ .
\endaligned
$$
Therefore
$$v\ \ddiv_{\xi.v}(D) =  \ddiv_\xi(D)v+D(v) \ ,
$$
and, multiplying both sides by $v^{-1}$, we obtain formula
(\ref{conformalchange}), since $v$ is even.
\qed

\subsection{Properties of generators defined by berezinian volumes}
We shall now assume that there is an odd Poisson structure, $\pi$,
on $(M,{\mathcal A})$, with odd Poisson bracket $\lcf\ ,\ \rcf$. 
Let $\xi$ be a berezinian volume on
$(M,{\mathcal A})$. Following the general pattern of Section
\ref{1.3}, we define the operator
$\Delta^{\pi,\xi}:{\mathcal A}\to {\mathcal A}$ by
\begin{equation} \label{genop}
\Delta^{\pi,\xi} f= (-1)^{|f|}\frac 12 \ddiv_\xi X^{\pi}_f \ ,
\end{equation}
for any section $f$ of ${\mathcal A}$. It follows from Proposition 
\ref{divberez}
and Theorem \ref{generateur0} that the odd operator 
$\Delta^{\pi,\xi}$ is a generator of bracket
$\pi$.
Thus, given an odd Poisson bracket, to any 
berezinian volume there corresponds a generator 
of this bracket. We shall now study the effect on the generator
of a change of berezinian volume, and determine 
under which conditions the
generator corresponding to a berezinian volume is of square 0.

\medskip

It follows from \eqref{particular} and Proposition \ref{prime}
that, when $\xi$ is a berezinian volume
and  $v=e^{2w}$
an invertible, even section of  ${\mathcal A}$,
\begin{equation}\label{newdelta}
\Delta^{\pi,\xi.v} = \Delta^{\pi,\xi} + X^{\pi}_w \ .
\end{equation}
It also follows from Proposition
\ref{genprime} that, if $\xi$ is a berezinian volume
such that $(\Delta^{\pi,\xi})^{2} = 0$, and $v=e^{2w}$ is
an invertible, even section of ${\mathcal A}$, then,
\begin{equation} \label{preMC}
(\Delta^{\pi,\xi.v})^2 = X^{\pi}_{\Delta^{\pi,\xi}w+
\frac 12 \lcf w,w\  \rcf} \ .
\end{equation}

Moreover,
$$ e^{-w}\Delta^{\pi,\xi}(e^{w})
= \frac 12 e^{-w} \ddiv_\xi(e^{w} X^{\pi}_w)
=
\Delta^{\pi,\xi}w
+
\frac 12
e^{-w} \lcf w,e^{w}\rcf \ ,
$$
and therefore
\begin{equation} \label{qqq}
\Delta^{\pi,\xi}w
+
\frac 12
\lcf w,w\rcf = e^{-w}\Delta^{\pi,\xi}e^{w} \ .
\end{equation}

In the context of supermanifolds (usually infinite-dimensional), the
Maurer-Cartan equation,
\begin{equation} \label{mastereq}
\Delta^{\pi,\xi} w + \frac 12 \lcf w,w \rcf = 0\ ,
\end{equation}
is referred to as the {\it Batalin-Vilkovisky quantum master equation}. 
In the case of odd symplectic
supermanifolds, the results stated below can be found in
articles dealing with
the BV-quantization of gauge theories, starting with \cite{B},
\cite{BV}, 
followed by, among others, \cite{W}, 
\cite{kh0}, \cite{Sc}, \cite{kh}, \cite{ASZK}, or of string
theories
\cite{Z}. See also \cite{G} and \cite{St}.
In our treatment, the more general 
case of possibly degenerate odd Poisson structures
is included.

\begin{proposition} \label{master}
Let $\xi$ be a berezinian volume on $(M,{\mathcal A})$
such that $(\Delta^{\pi,\xi})^{2} = 0$, and let $v=e^{2w}$ be
an invertible, even section of ${\mathcal A}$.

{\rm (}{\it i}{\rm )} The following conditions are equivalent
\begin{itemize}
   \item
$\Delta^{\pi,\xi}(e^{w})=0$,
   \item 
$w$ is a solution of \eqref{mastereq}.
\end{itemize}

{\rm (}{\it ii}{\rm )} If this condition is satisfied, then  
$(\Delta^{\pi,\xi.v})^{2}=0$.
\end{proposition}

\noindent{\it Proof.}
These implications follow from Equations 
\eqref{qqq} and \eqref{newdelta}. 
\qed

Conversely $(\Delta^{\pi,\xi.v})^{2}=0$ implies that there exists an
odd Casimir section $C$ of square 0 such that
$\Delta^{\pi,\xi}(e^{w}) C=0$.
(A Casimir section is a section
of ${\mathcal A}$ such that $X^{\pi}_C=0$.) In fact, if  
$(\Delta^{\pi,\xi.v})^{2}=0$, then there exists a Casimir section $C$ 
such that $e^{-w}\Delta^{\pi,\xi}(e^{w}) = C$, 
or $\Delta^{\pi,\xi}(e^{w})=e^{w}~C$. Together with
$(\Delta^{\pi,\xi})^{2}=0$,
this condition implies that 
$\Delta^{\pi,\xi}(e^{w}) C=0$, whence also $C^2=0$.

If the odd Poisson bracket $\pi$ is nondegenerate, any Casimir
section is a constant, therefore even, and necessarily $C=0$. In this
case, ({\it i}) and ({\it ii}) in the proposition are equivalent.

In quantum field theory, $e ^{{\frac {i}{\hbar}}S}$ is the action, and
the condition $\Delta^{\pi,\xi}(e ^{{\frac {i}{\hbar}}S})$  
states that the action is closed with respect to the
differential $\Delta^{\pi,\xi}$.

\begin{proposition} \label{newgenerator}
If $v=e^{2w}$ is an invertible, even section of $\mathcal A$ such that
$w$ is a solution of
the Equation {\rm (\ref{mastereq})},
then, for any section $f$ of ${\mathcal A}$,
$$ \label{deltaprime}
\Delta^{\pi,\xi.v} f = e^{-w}\Delta^{\pi,\xi}(e^{w}f) \ .
$$
\end{proposition}
\noindent{\it Proof.}
Using (\ref{BV}), we find that
$$
e^{-w}\Delta^{\pi,\xi}(e^{w}f) = e^{-w} \left( \lcf e^{w},f\rcf +
(\Delta^{\pi,\xi} e^{w})f +
e^{w}(\Delta^{\pi,\xi} f) \right) \ .
$$
By \eqref{qqq} and \eqref{newdelta},
$$ e^{-w}\Delta^{\pi,\xi}(e^{w}f)
= \Delta^{\pi,\xi.v}(f) + \left(
\Delta^{\pi,\xi}w
+
\frac 12 \lcf w,w\rcf \right) f \ .
$$
The result is proved in view of \eqref{newdelta}.
\qed

\medskip

\subsection{The supermanifold $\Pi T^{*}M$} \label{2.3}

For any manifold $M$ of dimension $n$, we consider
the supermanifold $\Pi T^{*}M$
of dimension $n|n$, whose structural sheaf is the sheaf of
multivectors on $M$.
The supermanifold $\Pi T^{*}M$
has an odd Poisson
bracket, the Schouten bracket of multivectors on $M$.
It is in fact nondegenerate, {\it i.e.}, the odd Poisson structure on
$\Pi T^{*}M$ is symplectic.
Here we revert to the usual notations for
vector fields, differential forms and functions on the ordinary manifold
$M$. See the Appendix for the definition to be used below of the
map $\alpha \mapsto \tilde{\alpha}$ from forms on a supermanifold
$(M,{\mathcal A})$ to
forms on $M$.

\begin{lemma}Given a volume form $\mu$ on $M$, there is a unique
berezinian volume $\xi_\mu$ on $\Pi T^{*}M$
such that
$$\widetilde{\xi_{\mu}(X)} = (i_{X_{(n)}}\mu)\mu \ , $$
for any field of multivectors $X$, where $i_{X_{(n)}}\mu$
is the result of
the duality-pairing of the homogeneous component $X_{(n)}$ of degree $n$
of the multivector $X$ and the $n$-form~$\mu$.
\end{lemma}

\noindent {\it Proof.} To any graded vector field, $D$, on $\Pi
T^{*}M$ is associated a
vector field, $\widetilde{D}$, on $M$, defined by
$$\widetilde{D}(f) = \widetilde{D(f)}$$ for any function $f$ on $M$.
Given a differential $n$-form, $\mu$, on $M$, we can define a graded 
$n$-form,
$\mu^G$, on $\Pi T^{*}M$ by $$<D_1,\dots,D_n ~ , ~ \mu^G>
=\mu(\widetilde{D_1},\dots,\widetilde{D_n}) \ ,$$ for graded vector
fields $D_1,\dots,D_n$. Then
$\widetilde{\mu^G} = \mu$.
Since the map $X \mapsto i_{X_{(n)}}\mu$
is a differential operator of order $n$ on the structural sheaf of
$\Pi T^{*}M$, the map
$\xi_{\mu} : X \mapsto (i_{X_{(n)}}\mu)\mu^G$
defines a section of the berezinian
sheaf, which is a berezinian volume if and only if $\mu$ is a volume
form. \qed

We remark that, for any positive function $v$ on $M$,
$$\xi_{v\mu} = (\xi_\mu).v^2.$$

\medskip

We assume that $M$ is an orientable
manifold, and we let $\mu$ be a volume form on $M$.
In the non orientable case, densities must be used instead of volume
forms.
On $\Pi T^{*}M$, there is an operator  $\Delta^{Schouten,\xi_{\mu}}$
associated to the odd
Poisson bracket and to the berezinian volume
$\xi_{\mu}$ by means of \eqref{genop}.
Let $\d$ be the de~Rham differential on $M$, and let $*_{\mu}$ be the
isomorphism from multivectors to forms defined by the volume form $\mu$.
Then we know (see, {\it e.g.,} \cite{KS}) that
$\partial_{\mu}= - *^{-1}_{\mu} \d *_{\mu}$ is a generator of the
Schouten bracket.
\bigskip

\begin{theorem}
\label{aa} 
For any volume form $\mu$ on $M$, 
the generator $\Delta^{Schouten,\xi_{\mu}}$ of the Schouten bracket
coincides with 
$\partial_{\mu} = - *_{\mu}^{-1} \d ~ *_{\mu} \ ,$ 
and $(\Delta^{Schouten,\xi_{\mu}})^{2}=0$.
\end{theorem}
\medskip

\noindent {\it Proof.} It is enough to show
that these operators coincide on
vector fields.
We prove this fact using local coordinates $(x^{1},...,x^{n},\xi_{1},
...,\xi_{n})$ on $\Pi T^{*}M$.
For a vector field $X=\sum_{i=1}^n X^{i} \frac {\partial}{\partial x^{i}}$
considered as a function $\sum_{i=1}^n X^{i}\xi_{i}$ on $\Pi T^{*}M$,
$$
\lcf X, ~ . ~ \rcf = \sum_{i=1}^{n} X^{i} \frac {\partial}{\partial
x^{i}} -
\sum_{i=1}^{n}\sum_{j=1}^{n}\frac {\partial X^{i}} {\partial x^{j}}
\xi_{i}
\frac {\partial} {\partial \xi_{j}} \ .
$$
Assume that $\mu_{0}= dx^{1}\wedge \ldots \wedge dx^{n}$,
so that
$$
- \frac 12 \ddiv_{\xi_{\mu_{0}} } \lcf X, ~ . ~ \rcf = - \frac 12
\sum_{i=1}^{n}
\left( \frac {\partial X^{i}}
{\partial x^{i}} +
\frac {\partial X^{i}}{\partial x^{i}}\right)= - \ddiv_{\mu_{0}}X =
\partial_{\mu_{0}}X \ .
$$
More generally, if $\mu = e^{w} \mu_{0}$,
then $\xi_{\mu} = \xi_{\mu_{0}}.e^{2w}$, and
for any vector field $X$,
$$
- \frac 12 \ddiv_{\xi_{\mu}} \lcf X, ~ . ~ \rcf
= - \frac 12 \ddiv_{\mu_{0}}X -X.w
=\partial_{\mu_{0}} X - X.w = \partial_{e^{w}\mu_{0}}X =
\partial_{\mu}X \ .
$$
The fact that $(\Delta ^{Schouten, \xi_{\mu}})^{2}=0$ 
is now an immediate consequence of the
fact that $\d ^{2}=0$. \qed

\medskip

\noindent{\it Remark.}
The equality $\Delta = \partial_{\mu}$ means that, for any differential
form
$\alpha$,
$$*^{-1}_{\mu} d \alpha = - \Delta *^{-1}_{\mu}\alpha \ . $$
The map $*^{-1}_{\mu}$
coincides with the ``Fourier transform with respect to the odd variables''
introduced in \cite{Vor} and \cite{Sc}, p. 255. Therefore,
in the case of a nondegenerate Poisson structure, our result
reduces to that of
\cite{Vor} and \cite{Sc}, formula (20), a fact
already observed by Witten
in \cite{W}, formula (13).
In the terminology of Voronov and Schwarz,
the operator $\partial_{\mu}$ on the functions
on $\Pi T^{*}M$
is the ``Fourier transform''
of the de~Rham differential acting on functions on $\Pi TM$.
Schwarz proves that a supermanifold $(M,{\mathcal A})$
of dimension $n|n$ with an odd
symplectic structure is equivalent in a suitable sense
to $\Pi T^{*}M$ with its canonical,
odd symplectic structure.

Now, let $\mu$ be a volume form, and $\omega$ a differential
form on $M$ such that $\omega_{(n)}$ is a volume form. If
$Q = *^{-1}_{\mu} \omega $,
then
$\xi' =\xi_{\mu}.Q^{2}$ is a berezinian volume on $\Pi T^{*}M$.
Setting $\Delta' = \Delta^{Schouten,\xi'}$, we see from Proposition
\ref{master} that $(\Delta')^{2}=0$
if $\Delta Q = 0$.
By Theorem \ref{aa}, this condition is equivalent to $\d \omega =0$. 
So
$(\Delta')^{2}= 0$ if $\omega$ is a closed form.
This result constitutes part of Theorem 5 of \cite{Sc}.
Moreover, it is proved there that
two closed forms in the same de~Rham cohomology class yield equivalent
structures.

\subsection{The supermanifold $\Pi TM$} \label{2.4}

We shall now consider another supermanifold 
attached to a smooth manifold $M$ of dimension $n$.
Let
${\mathcal A}=\Omega(M)$ be the sheaf of
differential forms on $M$. The pair $(M,{\mathcal A})$, usually 
denoted
by $\Pi TM$, is a supermanifold of dimension $n|n$.
The sections of the sheaf $\Omega(M)$ will be denoted by
$\alpha,\beta\ldots$, but differential
forms of ${{\mathbb Z}}$-degree $0$, {\it i.e.}, smooth functions,
will be denoted by, $f,g,\ldots$.
If $\alpha$ is a section of
$\Omega(M)$, then we denote
the homogeneous component of $\alpha$ of degree $n$ by $\alpha_{(n)}$.

\subsubsection{Canonical berezinian volume on the supermanifold $\Pi TM$}
\begin{lemma}
There is a unique berezinian volume $\xi$ on $\Pi TM$,
such that, for any section $\alpha$ of $\Omega(M)$,
\begin{equation}\label{ber}
\widetilde{\xi(\alpha)} = \alpha_{(n)} \ .
\end{equation}
 \end{lemma}

\noindent {\it Proof.}
If $\xi$ and $\xi'$ are berezinian volumes satisfying
\eqref{ber}, then
$\widetilde{(\xi-\xi')(\alpha)}=0$ for any section
$\alpha$ of $\Omega(M)$, and this
means, by the definition of the berezinian sheaf, that $\xi = \xi'$,
which proves the uniqueness of a berezinian volume satisfying
\eqref{ber}.
To prove its existence, we use local
coordinates and we show the invariance under a change of coordinates.
Let $(x^{1},\dots,x^{n})$ be local coordinates on an open set $U$ of the
manifold $M$. Then
$(x^{1},\dots,x^{n}, s^{1}=\d x^{1},\dots,s^{n}=\d x^{n})$
are graded local coordinates in $\Pi TM$, and a
local basis of derivations is
$$({\mathcal L}_{\frac\partial{\partial x^{1}}},\dots,
{\mathcal L}_{\frac\partial{\partial x^{n}}},
\frac\partial{\partial s^{1}} = i_{\frac\partial{\partial x^{1}}},\dots,
\frac\partial{\partial s^{n}} =i_{\frac\partial{\partial x^{n}}}) \ .$$
We now consider the local section of the berezinian sheaf,
$$
\xi_{U}= [\d^Gx^{1}\wedge\dots\wedge\d^Gx^{n} \otimes
i_{\frac\partial{\partial x^{1}}\wedge\dots\wedge
\frac\partial{\partial x^{n}}}] \ .
$$
A change of coordinates from $(x^{1},\dots,x^{n})$ to
$(y^{1},\dots,y^{n})$
induces a change of graded coordinates to
$(y^{1},\dots,y^{n}, t^{1}=\d y^{1},\dots,t^{n}=\d y^{n})$,
with matrix
$$
\begin{pmatrix}
\big(\frac{\partial y^{i}}{\partial x^{k}}\big) & 0\\
0 & \big(\frac{\partial y^{i}}{\partial x^{k}}\big)
\end{pmatrix} \ ,
$$
whose berezinian is equal to $1$.
Therefore, we can define a berezinian section $\xi$
on $\Pi TM$ by piecing the
locally defined $\xi_{U}$'s together.

This berezinian section is also a berezinian volume. In fact, this is true
locally and, if $U$ and $V$
are open sets such that $U \cap V \neq \emptyset$,
and if
local forms $\alpha$ and $\beta$ satisfy
$\xi_{U}.\alpha =\xi_{V}.\beta$, then $\alpha = \beta$, on $U\cap V$. 

Finally, if
$\alpha_{(n)} = f \d x^{1}\wedge\dots\wedge \d x^{n}$, where $f\in
C^\infty(U)$, then $\xi_{U}(\alpha) = f
\d^G x^{1}\wedge\dots\wedge \d^G x^{n}$, and therefore
$$\widetilde{\xi_{U}(\alpha)} =
\widetilde{f \d^G x^{1}\wedge\dots\wedge \d^G x^{n}} =
f \d x^{1}\wedge\dots\wedge \d x^{n} =
\alpha_{(n)} \ ,
$$
and $\xi$ satisfies condition \eqref{ber}.
\qed

\subsubsection{The canonical divergence operator on $\Pi TM$} \label{FNth}
Let us denote by $\ddiv_{can}$ the divergence operator associated to the
canonical berezinian volume on $\Pi TM$. The graded vector fields on
$\Pi TM$ are the derivations of the sheaf of differential forms.
Let us denote
the sheaf of vector-valued differential $k$-forms
by $\Omega^{k}(M;TM)$, for $k \geq 0$.
By the Fr{\"o}licher-Nijenhuis theorem \cite{FN}, we know that a derivation
$D$ of degree $k$ of $\Omega(M)$ can be
uniquely written as
$$D= {\mathcal L}_K+i_L \ ,$$
where $K$ is a section of $\Omega^k(M;TM)$ and $L$ is a section of
$\Omega^{k+1}(M;TM)$.

We introduce the notation ${\mathcal C}$ for the $(1,1)$-contraction
map from $\Omega^k(M;TM)$ to $\Omega^{k-1}(M)$, defined on a
decomposable element $K = \omega\otimes X$, where $X$ is a vector field
and $\omega$ is a $k$-form,
by ${\mathcal C}K = i_X\omega$, for $k \geq 1$, and by
${\mathcal C} =0$ on $\Omega^{0}(M;TM)$.

\begin{lemma}
For $K$ a section of $\Omega^k(M;TM)$ and
$L$ a section of $\Omega^{k+1}(M;TM)$,
$$\ddiv_{can}(i_{L}) = (-1)^k {\mathcal C}L \ ,\quad\ddiv_{can}(\d) =0\ ,
\quad\ddiv_{can}({\mathcal L}_K) = -\d({\mathcal C}K) \ ,$$ where $\d$ 
denotes the de~Rham
differential.
\end{lemma}

\noindent{\it Proof.}
We shall first compute $\ddiv_{can}(i_L)$ for a decomposable $L =
\omega\otimes X$ where $\omega$ is a section of $\Omega^{k+1}(M)$.
For any differential form $\alpha$,
$$\aligned
\widetilde{({\mathcal L}_{i_L}\xi)\alpha} &=
\widetilde{(\xi\circ i_L)\alpha}
 = \widetilde{\xi(i_L\alpha)} \\
&= \widetilde{\xi(\omega\wedge i_X\alpha)} =
(\omega\wedge i_X\alpha)_{(n)}=
\omega\wedge i_X(\alpha_{(n-k)})\\
& = (-1)^k i_X(\omega)\wedge\alpha_{(n-k)}
= (-1)^k (i_X(\omega)\wedge\alpha)_{(n)}\\
&= \widetilde{\xi((-1)^k i_X(\omega)\wedge\alpha)} \ ,
\endaligned$$
where we have used the relation
$$(i_X\omega)\wedge\alpha_{(n-k)} + (-1)^{k+1}\omega\wedge
i_X\alpha_{(n-k)} = i_X(\omega\wedge\alpha_{(n-k)}) = 0 \ .$$
Therefore ${\mathcal L}_{i_L}\xi = \xi.((-1)^k i_X(\omega)) =
\xi.((-1)^k {\mathcal C}L)$, and $\ddiv_{can}(i_L) = (-1)^k {\mathcal 
C}L$.

Similarly, for the derivation $\d$,
$$
\widetilde{({\mathcal L}_{\d}\xi)\alpha} =
\widetilde{(\xi\circ \d)\alpha} =
\widetilde{\xi(\d\alpha)} \\
= \d(\alpha_{(n-1)}) \ ,$$
which is always an exact form. Thus ${\mathcal L}_{\d}\xi = 0$.
Finally, using \eqref{commutator}, we obtain
$$\ddiv_\xi({\mathcal L}_K) = \ddiv_\xi([i_K,\d]) = 
i_K(\ddiv_\xi(\d)) -
(-1)^{k-1}\d(\ddiv(i_K)) = -\d({\mathcal C}K) \ .\qed$$

In particular, if $X$ is a vector field on $M$, then
\begin{equation}\label{hh}
\ddiv_{can}(i_X) = 0 \ , \quad \ddiv_{can}({\mathcal L}_X) = 0 \ .
\end{equation}

\subsubsection{
Generators of the Koszul-Schouten bracket}
We shall now assume that the base manifold $M$ is equipped with 
a Poisson structure.
Given a Poisson manifold $(M,P)$, there is an odd Poisson bracket,
$\lcf\ ,\ \rcf_P$, on the supermanifold $\Pi TM$, called the
{\it Koszul-Schouten bracket}, that is characterized by the conditions,
$$\lcf f,g\rcf_P = 0 \ , \qquad
\lcf f, \d g \rcf_P = \{ f,g\} \ , \qquad
\lcf \d f, \d g \rcf_P = \d \{f,g\} \ ,
$$
for all $f,g\in C^\infty(M)$, where $\{ ~ , ~ \}$ denotes the Poisson
bracket on $C^\infty(M)$ defined by $P$, together with
the graded Leibniz rule.  It was shown by Koszul \cite{Kz1} that
a generator for this bracket is the {\it Poisson homology
operator},
$\partial_P = [\d,i_P]$,
sometimes called the Koszul-Brylinski operator. See \cite{Hu1}, and 
also \cite{KS}.
On the other hand, we know that, given a berezinian volume $\xi$ on $\Pi
TM$,  the operator $\Delta^{(P),\xi}$ defined by
$$\Delta^{(P),\xi}(\alpha) = (-1)^{|\alpha|} \frac 12
\ddiv_\xi(\lcf\alpha,\ . \ \rcf_P) \ ,$$
for any differential form $\alpha$,
is also a generator of $\lcf\ ,\ \rcf_P$.

\begin{theorem} \label{bb}
The generator $\Delta^{(P),can}$
of the Koszul-Schouten bracket associated to
the canonical berezinian volume coincides with
$\partial_P$, and $(\Delta^{(P),can})^2=0$.
\end{theorem}
\noindent{\it Proof.}
It suffices to prove that $\Delta^{(P),can}$ and $\partial_P$ agree on
$1$-forms, and it is enough to show that both vanish on exact
$1$-forms, $\alpha = \d f$, where $f \in C^{\infty}(M)$.
In fact,
$\lcf \d f,\ . \ \rcf_P = {\mathcal L}_{\#_P\d f}$,
where $ \#_P\d f = \{f,~.~\} \ .$
From \eqref{hh}, it follows that
$\Delta^{(P),can}(\d f) =
 - \frac 12 \ddiv_{can}({\mathcal L}_{\#_P\d f}) = 0 \ .$ 
And clearly $\partial_P(\d f) = 0 \ .$
Moreover, $(\partial _P)^2=0$, since $[\d ,\partial _P]= 0$ and
$[i_P,\partial _P]=0$, and therefore $(\Delta^{(P),can})^2=0$.
\qed

\medskip

\noindent{\it Remark.}
Any nondegenerate metric $g$ on the manifold $M$ defines an isomorphism
from multivectors to differential forms. Hence,
from the Schouten bracket of multivectors,
we obtain a $\mathbb Z$-graded bracket $\lcf~,~\rcf_{g}$ on the
sheaf of differential
forms on $M$.
Then, the codifferential $\delta_{g}$ associated to $g$ is a generator of 
this
bracket. (See \cite{V94} or \cite{BMS-V}.) One can also consider
the operator associated to the canonical berezinian on $\Pi TM$, defined 
by
$$\Delta^{(g),can}(\alpha) = (-1)^{|\alpha|} \frac 12
\ddiv_{can}(\lcf\alpha,~ . ~ \rcf_g) \ ,$$
for any differential form $\alpha$,
and one can show that these two generators
of the bracket $\lcf\ ,\ \rcf_g$ coincide,
$\Delta^{(g),can} = \delta_{g}$.

\subsection{QS, SP and QSP-manifolds} \label{2.5}

The following definitions, adapted from \cite{Sc} and \cite{Sc2},
will be useful in order to reformulate some of our results.

\begin{definition}
Let $(M,{\mathcal A})$ be a supermanifold, $D$ an odd vector field
and $\xi$ a berezinian volume. We say that
$((M,{\mathcal A}),D,\xi)$ is a {\rm QS-manifold} if
$D^2 = 0$ and $\ddiv_{\xi}(D) = 0$.
\end{definition}

\begin{definition}
Let $(M,{\mathcal A})$ be a supermanifold, 
$\pi$
an odd Poisson bracket and $\xi$ a berezinian volume. 
We say that
$((M,{\mathcal A}),\pi, \xi)$ is a {\rm weak SP-manifold} if
$(\Delta^{\pi,\xi})^2 = 0$, where $\Delta^{\pi,\xi}$ is defined by
\eqref{genop}.
If the Poisson bracket is nondegenerate, then the supermanifold is
an {\rm SP-manifold}.
\end{definition}

\begin{definition}
Let $(M,{\mathcal A})$ be a supermanifold, $\pi$ an odd Poisson
bracket, $D$ an odd vector field and $\xi$ a berezinian volume. We
say that
$((M,{\mathcal A}),\pi,D,\xi)$ is a {\rm weak QSP-manifold} if
\begin{itemize}
\item $((M,{\mathcal A}),\pi,\xi)$ is a weak SP-manifold,
\item $((M,{\mathcal A}),D,\xi)$ is a QS-manifold, and
\item $D$ is a derivation of the odd Poisson bracket, $\pi$.
\end{itemize}
If the Poisson bracket is nondegenerate, and if $D$ is the
hamiltonian vector field defined by an even section of $\mathcal A$, 
then the supermanifold is
a {\rm QSP-manifold}.
\end{definition}

If the graded vector field $D$ is the hamiltonian vector field 
$\lcf h, ~. ~ \rcf $, where $h$ is an even section of
${\mathcal A}$, then $D^2=0$ if and only if $\lcf h,h \rcf= 0 \ .$
In field theory, 
this condition appears under the name 
{\it classical master equation}.

It follows from identity \eqref{derdelta} that, 
if $\Delta$ is a generator of the odd
Poisson bracket $\lcf ~,~\rcf$ and if $D$ is a graded vector field on
$(M,{\mathcal A})$, a necessary and
sufficient condition for $D$ to be a derivation of the odd Poisson
bracket is that the graded commutator, $[D, \Delta]$,
be a derivation of
the associative multiplication of ${\mathcal A}$. This implies
\begin{proposition}
If $\pi$ and $\xi$ define a weak SP-structure
on $(M,{\mathcal A})$, if
$D$ and $\xi$ define a
QS-structure on $(M,{\mathcal A})$, and if $[D,\Delta]=0$, then
$\pi$, D and $\xi$ define a weak QSP-structure on $(M,{\mathcal A})$.
\end{proposition}

The following theorems follow in part from the 
results of Section \ref{2.3} and \ref{2.4}. See also \cite {ASZK}.

\begin{theorem} \label{piT*M}
({\it i}) For any manifold $M$ with a volume element,
the supermanifold
$\Pi T^{*}M$, with the Schouten bracket and the berezinian volume
$\xi_{{\mu}}$,
is an SP-manifold.

({\it ii}) Let $(M,P)$ be a Poisson manifold, 
and let $d_P = \lcf P, \ . \ \rcf$ be the
Lichnerowicz-Poisson
differential. 
Then $\Pi T^{*}M$, with the Schouten bracket, the odd vector field $d_P$, 
and the canonical berezinian volume, is a QSP-manifold. 
\end{theorem}

\noindent{\it Proof.}
In fact, the odd vector field $d_P$ is of square $0$, because
$\lcf P,P \rcf =0$, and it is a
derivation of the Schouten bracket by the graded Jacobi identity.
\qed

\begin{theorem} \label{piTM}
({\it i}) For any manifold $M$, 
the supermanifold $\Pi TM$, with the de~Rham differential $\d$ and the 
canonical berezinian volume,
is a QS-manifold.

({\it ii}) Let $(M,P)$ be a Poisson manifold, 
and let $\lcf\ ,\ \rcf_P$ be the
Koszul-Schouten bracket. Then $\Pi TM$ with 
the odd Poisson bracket $\lcf\ ,\ \rcf_P$, the de~Rham
differential $\d$ and the canonical berezinian 
volume, is a weak QSP-manifold.
If $P$ is a nondegenerate Poisson structure, then $\Pi TM$ is a
QSP-manifold.
\end{theorem}

\noindent{\it Proof.}
To complete the proof, we only have to recall that the
de~Rham differential $\d$ is a derivation of  $\lcf\ ,\ \rcf_P \ ,$
and that, when $P$ is nondegenerate
with inverse the symplectic form, $\omega$, then
$\d = \lcf \omega, ~ . ~ \rcf_{P} \ ,$
so that $\d$
is the hamiltonian vector field associated to
$\omega$.
(See \cite{KSM} or \cite{BM} for a proof of this fact.)\qed

\section{Linear connections and generators of odd Poisson
brackets} \label{3}

Divergence operators on smooth
manifolds
can be defined not only by means of
volume forms, but also by means of connections. (See \cite{Hu}.)
While in Section \ref{2}, we generalized the first approach to supermanifolds,
replacing
volume forms
by their graded analogue, the berezinian volumes, in this section we
generalize the second method, defining divergence operators by means of
graded connections.

\subsection{Divergence operators defined by graded connections} \label{divop}
We first recall the notion of graded
connection. See \cite{Ma} and \cite{MP} for the definitions of left and 
right graded connections.
Here we consider only left graded connections, which we simply
call connections. Let $(M,{\mathcal A})$ be a
supermanifold and let $\Der{\mathcal A}$ be the sheaf of derivations of 
${\mathcal A}$.

\begin{definition}
Let ${\mathcal S}$ be a sheaf of ${\mathcal A}$-modules on $M$.
A {\em left graded connection}, or simply a {\em connection},
on ${\mathcal S}$ is a morphism of
sheaves of graded vector spaces
from $\Der{\mathcal A} \otimes {\mathcal S}$ to ${\mathcal S}$, denoted
$D \otimes \alpha \mapsto \nnabla_{D}\alpha \ , $
which satisfies the identity
$$
\nnabla_{fD} \alpha = f \nnabla_{D}\alpha \ ,
$$
and the Leibniz rule,
$$
\nnabla_D (f\alpha) = D(f) \alpha + (-1)^{|D||f|} f\nnabla_D \alpha \ ,
$$
for any section $f$ of ${\mathcal A}$, any derivation $D$ of ${\mathcal 
A}$,
and any section $\alpha$ of ${\mathcal S}$.

A connection on the sheaf $\Der{\mathcal A}$ of derivations of $\mathcal A$
is a {\em graded linear connection} or simply a {\em linear
connection} on $(M,{\mathcal A})$.
\end{definition}

\begin{definition}
The {\em curvature}, $R^{\nnabla}$, of a connection,
$\nnabla$, on ${\mathcal S}$ is defined by
$$R^{\nnabla}(D_1,D_2)= [\nnabla_{D_1},\nnabla_{D_2} ] -
\nnabla_{[D_1,D_2]} \ ,$$
for any derivations $D_{1},D_{2}$ of ${\mathcal A}$.

The {\em torsion}, $T^{\nnabla}$, of a linear connection, $\nnabla$,
on $(M,{\mathcal A})$
is defined by
$$T^{\nnabla}(D_1,D_2)
= \nnabla_{D_1}D_2 - (-1)^{|D_1||D_2|}\nnabla_{D_2}D_1 -
[D_1,D_2] \ .$$
\end{definition}

We shall now define the divergence operator associated with a
linear connection on $(M,{\mathcal A})$. 
Let $sTr$ denote the supertrace
of an endomorphism of sheaves of ${\mathcal A}$-modules
(see, {\it e.g.}, \cite{Ma} or \cite{D}),
and let ${\rm ad}_{D}$ denote the endomorphism of $\Der{\mathcal A}$,
$E \mapsto [D,E]$.
For any graded vector field $D$, we set
\begin{equation}\label{divconn}
\ddiv_\nnabla(D) =  sTr\left(\nnabla_D - {\rm ad}_{D}\right) \ .
\end{equation}

\begin{proposition} \label{propdivconn}
For any linear connection, $\nnabla$,
on $(M,{\mathcal A})$, the map, 
$\ddiv_\nnabla:\Der {\mathcal A} \to {\mathcal A},$
defined by \eqref{divconn} is a 
divergence operator.
\end{proposition}

\noindent{\it Proof.}
The map $\ddiv_{\nnabla}$ is even. It follows from
$[f D, E] = f[D,E] -(-1)^{(|f|+|D|)|E|} E(f)D$, that
$$\aligned
\ddiv_\nnabla(f D)& =
 sTr(\nnabla_{f D} - {\rm ad}_{fD})\\
&=
f\ \ddiv_\nnabla(D)
+ sTr\left(E \mapsto (-1)^{(|f|+|D|)|E|} E(f)D\right) \ . \\
\endaligned
$$
Let the graded dimension of the supermanifold be $m|n$,
and let us choose a system of local graded coordinates
$(x^{1},\dots,x^{m},s^{1},\dots,s^{n})$. We find that
$$\aligned 
sTr(E &\mapsto (-1)^{(|f|+|D|)|E|} E(f)D)
=
\sum_{i=1}^{m}{\frac {\partial f}{\partial x^{i}}} D(x^{i})
- (-1)^{|f|+|D|} \sum_{\rho=1}^{n}
{\frac {\partial f}{\partial s^{\rho}}} D(s^{\rho})
\\&=
(-1)^{|f||D|} \sum_{i=1}^m  D(x^{i}){\frac {\partial f}{\partial
{x^{i}}}}
+  (-1)^{|f||D|}\sum_{\rho=1}^n  D(s^{\rho})
{\frac {\partial f} {\partial {s^{\rho}}}}
= (-1)^{|f||D|} D(f) \ ,
\endaligned
$$
where we have used the local expression of the derivation $D$ in the
basis
of local graded vector fields,
$\left(
{\frac {\partial}{\partial x^{1}}},
\dots,
{\frac {\partial}{\partial x^{m}}},
{\frac {\partial}{\partial s^{1}}},
\dots,
{\frac
{\partial}
{\partial s^{n}}} \right)$.
\qed

\begin{proposition} \label{curvtr}
Let $\nnabla$ be a torsionless
linear connection on ${\mathcal A}$
and let $D_1$ and $D_2$ be graded vector fields. Then
\begin{equation} \label{dd}
{\mathcal R}^{\ddiv_\nnabla} (D_1,D_2)= - sTr(R^\nnabla(D_1,D_2)) \ .
\end{equation}
\end{proposition}

\noindent{\it Proof.}
We have to prove that
$$\begin{aligned} 
D_1(sTr(\nnabla_{D_2}- 
{\rm ad}_{D_2}))&-(-1)^{|D_1||D_2|}D_2
(sTr(\nnabla_{D_1}-
{\rm ad}_{D_1}))-
sTr(\nnabla_{[D_1,D_2]}- {\rm ad}_{[D_1,D_2]}) \\ \nonumber
&= sTr([\nnabla_{D_1},\nnabla_{D_2}] - \nnabla_{[D_1,D_2]}).
\end{aligned}
$$
This result follows from a computation of these two expressions
in local coordinates, for pairs of
commuting graded
vector fields, $D_1$ and $D_2$,
in a local basis.
\qed

\subsection{Generators defined by graded connections} \label{3.2}
We now assume that $(M,{\mathcal A})$ has an odd Poisson structure, $\pi$,
whose
bracket we denote by $\lcf ~,~ \rcf$.
Let $\nnabla$ be a linear connection on $(M,{\mathcal A})$. Following the
general pattern of Section \ref{1.3}, we define the 
operator
$\Delta^{\pi,\nnabla}
:{\mathcal A}\to {\mathcal A}$ by
\begin{equation} \label{ppp}
\Delta^{\pi, \nnabla}(f) = (-1)^{|f|}\frac 12 \ddiv_\nnabla(\lcf f,\ .
\  \rcf) \ ,
\end{equation}
for any section $f$ of ${\mathcal A}$.
It follows from Proposition \ref{propdivconn} and Theorem \ref{generateur0}
that the odd operator $\Delta^{\pi,\nnabla}$ is a generator of bracket
$\pi$.
Therefore, to any linear connection on an odd Poisson manifold,
there corresponds a generator of the odd Poisson bracket.

\begin{proposition} \label{jj}
Let $\nnabla$ be a torsionless linear connection on $\mathcal A$. The
following properties are equivalent
\begin{itemize}
 \item
$\Delta^{\pi,\nnabla}$ is a derivation of the odd Poisson
bracket,$\pi$,
 \item
$(\Delta^{\pi, \nnabla})^{2}$ is a derivation of the sheaf of
associative algebras, ${\mathcal A}$,
 \item
$sTr(R^\nnabla)$ vanishes on the sheaf of hamiltonian derivations.
\end{itemize}
\end{proposition}

\noindent{\it Proof.}
These equivalences follow from Lemma \ref{ccc}, and from Corollary 
\ref{dercurv}
together with Proposition \ref{curvtr}.
\qed

We now compare the generators associated to torsionless
linear connections,
$\nnabla$ and $\nnabla'$, on $\mathcal A$. The difference
$\nnabla'_{D} - \nnabla_{D}$ is then a morphism of sheaves of
$\mathcal A$-modules from $\Der {\mathcal A}$ to itself, which we denote by
$u(D)$.

\begin{proposition}
Let $\nnabla$ and $\nnabla'$ be torsionless linear
connections on $\mathcal A$. Then
$$\Delta^{\pi, \nnabla'}(f) = \Delta^{\pi, \nnabla}(f)
+ (-1)^{|f|}\frac 12 sTr(u (X^{\pi}_{f})) \ , $$
for any section $f$ of ${\mathcal A}$ \ .
\end{proposition}
\noindent{\it Proof.}
This relation follows from the fact that, for any derivation $D$ of
$\mathcal A$, $\ddiv_{\nnabla'}(D) = \ddiv_{\nnabla}(D) + sTr(u(D)).$\qed

\medskip

\noindent{\it Remark.} 
In the case of an ordinary manifold, the trace of the curvature of a
linear connection is the curvature of the connection induced on the bundle of
top-degree forms. It would be interesting to interpret the
supertrace of the curvature of a graded linear connection as the curvature of
a connection on the sections of the berezinian sheaf.

\subsection{Metrics and metric connections on supermanifolds} \label{3.3}
We recall the definitions of metrics and metric linear
connections on supermanifolds.

\begin{definition} A {\rm graded metric}, or simply a {\rm metric},
on $(M,{\mathcal A})$
is a morphism of sheaves of
${\mathcal A}$-modules,
$\langle ~ , ~ \rangle:\Der\ {\mathcal A} \otimes
\Der\ {\mathcal A}\to {\mathcal A} $, such that
\begin{itemize}
\item $\langle D_1,D_2\rangle
= (-1)^{|D_1||D_2|}\langle D_2,D_1\rangle \ ,$
for derivations $D_{1}$ and $D_{2}$ \quad (graded symmetry),
\item the map $D\mapsto \langle D, ~ \cdot ~ \rangle$ is an isomorphism of
sheaves of ${\mathcal A}$-modules from $\Der {\mathcal A}$ to
${\rm Hom}_{\mathcal A}(\Der{\mathcal A},{\mathcal A})$ \quad (nondegeneracy).
\end{itemize}
\end{definition}

\begin{definition}
A linear
connection $\nnabla$ on ${\mathcal A}$ 
is {\rm metric} with respect to a metric 
$\langle~,~\rangle$ if, for any
derivations 
$D,D_1$, and $D_2$ of ${\mathcal A}$,
$$D\langle D_1,D_2\rangle = \langle\nnabla_DD_1,D_2\rangle +
(-1)^{|D_1||D|}\langle D_1,\nnabla^0_DD_2\rangle
+(-1)^{|D_1|(|D|+1)}\langle D_1,\nnabla^1_DD_2\rangle,$$
where $\nnabla = \nnabla^0+\nnabla^1$ is the decomposition of the linear
connection into its even and odd components.
 \end{definition}

The proof of the following theorem,
can be found in \cite{L0}, p. 134, and in 
\cite{MS}.

\begin{theorem} \label{LC}
There exists a unique torsionless
linear connection which is metric with respect to a given metric.
It is determined by 
$$\begin{array} {rcl}  
2 \langle\nnabla_{D_1}D_2,D_3\rangle 
&=&
D_1\langle D_2,D_3\rangle
+\langle [D_1,D_2],D_3\rangle \\
&+& (-1)^{|D_1|(|D_2|+|D_3|)} \left( D_2\langle D_3,D_1\rangle
- \langle [D_2,D_3],D_1\rangle \right) \\
&-&  (-1)^{|D_3|(|D_1|+|D_2|)} \left( D_3\langle D_1,D_2\rangle
- \langle [D_3,D_1],D_2\rangle \right) \ .
\end{array}
$$
\end{theorem}
The linear connection defined in Theorem \ref{LC} is
called the graded Levi-Civita connection or simply, the
{\em Levi-Civita connection} of the metric
$\langle ~,~ \rangle$. 
The Levi-Civita connection of a homogeneous metric
is even.
(See also \cite{dW}, where an
expression in local coordinates of the
Levi-Civita connection is given in the case of a
homogeneous, even metric.)

\subsection{Linear connections and Schouten bracket} \label{3.4}

We shall again consider the
supermanifold $\Pi T^{*}M$, whose sheaf of functions is 
the sheaf of multivectors on $M$.
We shall use the notations $f,g\dots$ for functions on $M$ or on an open
set of $M$, and the notations $X, Y, \dots$ for vector fields and
$\alpha,\beta\dots$ for differential $1$-forms.

\subsubsection{Graded vector fields on $\Pi T^{*}M$}
Whereas the derivations of the algebra of forms on a manifold can be
classified by the Fr{\"o}licher-Nijenhuis theorem \cite{FN} (and
see Section \ref{FNth}), the classification of the derivations of the
algebra of multivectors on $M$ requires the use of an auxiliary linear
connection, $\nabla$.
Let $U$ be an open set of $M$.
If $K= Q \otimes X$ is a vector-valued multivector on $U$, where $Q$
is a multivector and $X$ is a vector, we define
$\nabla_K V = Q\wedge \nabla_XV$, for any multivector $V$ on
$U$.
If $L = W \otimes \alpha$
is a $1$-form-valued multivector on $U$, where $W$ is a multivector
and $\alpha$ is a differential $1$-form, we define
$i_L V = W\wedge i_\alpha V$.

\begin{proposition} \label{jjj} \cite{MM}
Let $D$ be a graded vector field of degree $r$
on $\Pi T^{*}M$, {\it i.e.},
a derivation of degree $r$ of the sheaf of multivectors on $M$.
On any open set $U$ of $M$, there exist a vector-valued $r$-vector, $K$,
and a $1$-form-valued $(r+1)$-vector, $L$, each uniquely defined, such
that
$D|_{U} = \nabla_K+i_L.$
\end{proposition}

As a consequence, we see that, if $(e_{1},\dots,e_{n})$ is
a local basis of vector fields on $U$
and $(\epsilon^{1},\dots,\epsilon^{n})$ is the
dual basis, then
$(\nabla_{e_{1}},\dots,\nabla_{e_{n}},
i_{\epsilon^{1}},\dots,i_{\epsilon^{n}})$
generate the derivations of the algebra of multivectors over $U$,
as a module over the algebra of multivectors
over $U$.

\subsubsection{The graded connection on $\Pi T^{*}M$ associated
to a linear connection on $M$}
We shall show how to associate a metric on $\Pi T^{*}M$ to a linear
connection on $M$, and we shall study the Levi-Civita connection of
this metric.
\begin{definition} \label{def-metr}
Let $\nabla$ be a linear connection on $M$.
We define a metric
$\langle~,~\rangle_\nabla$ on $\Pi T^{*}M$ by its value
on derivations of type $\nabla_X$, where $X$ is
a vector field, and of type $i_\alpha$,
where $\alpha$ is a $1$-form,
$$
\langle \nabla_X,\nabla_Y\rangle_\nabla = 0,\quad
\langle \nabla_X,i_\alpha\rangle_\nabla = \alpha(X),\quad
\langle i_\alpha,i_\beta\rangle_\nabla = 0 \ .
$$
\end{definition}

To verify the nondegeneracy of the metric thus defined,
we observe that,
in the local basis of derivations
$(\nabla_{e_{1}},\dots,\nabla_{e_{n}},
i_{\epsilon^{1}},\dots,i_{\epsilon^{n}})$,
the matrix of this metric is $\begin{pmatrix} 0& Id \\ Id &0  \end{pmatrix}$.
This metric is odd.

\begin{proposition}\label{LC-odd}
Let $\nabla$ be a torsionless linear connection on $M$.
The Levi-Civita connection, $\nnabla$, of 
the metric $\langle ~,~\rangle_\nabla$ on $\Pi T^{*}M$ is given by
$$
\nnabla_{\nabla_X}\nabla_Y = \nabla_{\nabla_XY}+i_{R(~\cdot ~ ,Y)X} \ ,
\quad \nnabla_{\nabla_X}i_\alpha = i_{\nabla_X\alpha} \ ,\quad
\nnabla_{i_\alpha} = 0 \ ,
$$
where $R$ denotes the curvature tensor of $\nabla$.
\end{proposition}

\noindent{\it Proof.}
We shall make use of the commutation relations
$$ [\nabla_X,\nabla_Y] = \nabla_{[X,Y]} + i_{R(X,Y)},\quad
[\nabla_X,i_\alpha] = i_{\nabla_X\alpha}, \text{ and   }
[i_\alpha,i_\beta] = 0 \ ,$$
Using Theorem
\ref{LC}, Definition
\ref{def-metr} and the fact that the connection $\nabla$
is torsionless, we obtain

$$
\begin{array}{rcl}
\langle\nnabla_{\nabla_X}\nabla_Y,\nabla_Z\rangle_{\nabla} &=&
{\frac 12} ( R(X,Y)Z -R(Y,Z)X + R(Z,X)Y)\\
&=&  R(Z,Y)X=  \langle i_{R(~.~ ,Y)X},\nabla_Z\rangle_{\nabla} \ , \\
&&\\
\langle\nnabla_{\nabla_X}\nabla_Y,i_\alpha\rangle_{\nabla} &=&
 \alpha(\nabla_XY)= \langle
\nabla_{\nabla_XY},i_\alpha\rangle_{\nabla} \ ,
\end{array}
$$
and
$$
\langle\nnabla_{\nabla_X}i_\alpha,\nabla_Y\rangle_{\nabla} =
(\nabla_X\alpha)Y=  \langle
i_{\nabla_X\alpha},\nabla_Y\rangle_{\nabla },\quad
\langle\nnabla_{\nabla_X}i_\alpha,i_\beta\rangle_{\nabla} = 0 \ .
$$
From these relations and the nondegeneracy of the graded metric we
obtain the first two formul{\ae}, while the third follows
from the fact that
$\nnabla$ is torsionless.\qed

\begin{proposition}
The curvature $R^{\nnabla}$ of 
the Levi-Civita connection, $\nnabla$, of
the metric
$\langle~,~\rangle_\nabla$ on $\Pi T^{*}M$, satisfies
$$
\aligned
R^\nnabla(\nabla_X,\nabla_Y)\nabla_Z &=
\nabla_{R(X,Y)Z}+i_{(\nabla_XR)(~\cdot ~ ,Z)Y}+i_{(\nabla_YR)(~\cdot
~,Z)X} \ ,\\
R^\nnabla(\nabla_X,\nabla_Y)i_\alpha &= i_{R(X,Y)^\ast\alpha} \ ,
\endaligned
$$
where $R$ denotes the curvature tensor of $\nabla$,
and $R(X,Y)^{*}$ denotes the transpose of
$R(X,Y)$. Moreover
$R^{\nnabla}(\nabla_{X},i_{\alpha})=R^{\nnabla}(i_{\alpha}, i_{\beta})=0$.
\end{proposition}

\noindent{\it Proof.}
The proof is a straightforward
computation using Proposition \ref{LC-odd}.\qed

\begin{corollary} \label{vvv}
Let $\nabla$ be a torsionless linear connection on $M$. Then the 
Levi-Civita
connection of the metric
$\langle~,~ \rangle_\nabla$ on $\Pi T^{*}M$
is flat
if and only if $\nabla$ is flat.
\end{corollary}

\subsubsection{Generators of the Schouten bracket} \label{gensch}

We have just seen that, to a
torsionless linear connection $\nabla$ on $M$
we can associate the
Levi-Civita connection, $\nnabla$, of the odd metric
$\langle~,~\rangle_\nabla$ on $\Pi T^{*}M$, and therefore, by 
\eqref{ppp},
a generator of the Schouten bracket, which we shall denote by 
$\Delta^{Schouten,\nnabla}$.
There exists another construction, due to Koszul \cite{Kz1},
which associates a generator of the Schouten bracket 
to a torsionless linear connection 
$\nabla$ on $M$. 
To $\nabla$,
he first associates the corresponding {\it divergence
operator}, defined on vector fields $X$ by 
\begin{equation} \label{trace}
\ddiv_{\nabla}(X)= Tr(\nabla_{X}-{\rm ad}_{X}) \ .
\end{equation}
This is the definition that
is used in fact in \cite{Kz1} (although it appears by mistake
with the opposite sign in
its first occurrence, page 262, before Lemma (2.1)).
This map is a divergence operator, {\it i. e.}, satisfies
\eqref{gooddiv},
on the purely even algebra
$C^{\infty}(M)$. 
For a flat connection on flat space, it reduces to the
elementary divergence. 
He then shows, using a local basis of vector fields, that 
there is a unique operator on the multivectors,
${\Delta}^{\nabla}$,
of degree $-1$, that extends the operator $- \ddiv_{\nabla}$
and generates the Schouten bracket.
We shall now show that the generators of the Schouten bracket obtained
by these two constructions coincide.

\begin{lemma} \label{OOO}
For any vector field $X$, and for any $1$-form $\alpha$,
$$\ddiv_\nnabla(\nabla_X) = \ddiv_\nabla(X),\quad
\ddiv_\nnabla(i_\alpha) = 0 \ .$$
\end{lemma}

\noindent{\it Proof.}
If $(x^{1},\dots,x^{n})$ is a system of local coordinates on $M$,
then a local basis of graded derivations on $\Pi T^\ast M$ is
$(\nabla_{\frac
\partial{\partial x^{1}}},\dots, \nabla_{\frac
\partial{\partial x^{n}}}, i_{\d x^{1}},\dots,i_{\d x^{n}}) \ .$
We use the relations
$<\nabla_{\frac \partial{\partial x^{j}}},\d ^G x^k> = \delta_{j}^{k} \ ,$
and  $<i_{\d x^{j}},\d^G x^k> = 0 \ .$
In order to compute
$\ddiv_{\nnabla}(\nabla_{X}) = sTr(\nnabla_{\nabla_X}-{\rm
ad}_{\nabla_X})$, we first observe that, because
$\nnabla$ is torsionless,
$$
\nnabla_{\nabla_X}i_{\d x^{j}}-[\nabla_X,i_{\d
x^{j}}] =  \nnabla_{i_{\d x^{j}}}\nabla_X = 0 \ .
$$
Therefore 
$$\begin{aligned}
&sTr(\nnabla_{\nabla_X} -{\rm ad}_{\nabla_X})
= \sum_{j= 1}^n < \nnabla_{\nabla_X}\nabla_{\frac
\partial{\partial x^{j}}} -[\nabla_X,\nabla_{\frac
\partial{\partial x^{j}}}] , \d^{G}x^{j} >\\
&= \sum_{j= 1}^n <\nnabla_{\nabla_{\frac
\partial{\partial x^{j}}}}\nabla_X , \d^{G}x^{j} >
= \sum_{j= 1}^n <\nabla_{\frac
\partial{\partial x^{j}}}X , \d x^{j} >
= \ddiv_\nabla(X) \ .
\end{aligned}
$$
\qed
\begin{theorem} \label{cc}
For any torsionless linear connection $\nabla$ on $M$,
the generators $\Delta^{Schouten,\nnabla}$  and 
${\Delta}^{\nabla}$ of the Schouten
bracket coincide.
If $\nabla$ is flat, this generator is of square $0$. 
\end{theorem}

\noindent{\it Proof.}
Since we know that both operators
are generators of the Schouten bracket,
we need only show that they coincide on functions and on vector
fields.  On functions, $\Delta^{Schouten,\nnabla}$ vanishes since
$\lcf f, ~\cdot~ \rcf = i_{df}$ and $\ddiv_\nnabla(i_{df}) = 0$, as
does
$\Delta^{\nabla}$ because it is of degree $-1$.
Now, for any vector field $X$ and any torsionless
linear connection $\nabla$ on
an open set $U$ of $M$,
$\lcf X, ~.~\rcf= \nabla_{X}-i_{\nabla X}$, since both derivations of the
sheaf of multivectors coincide on functions and on vectors.
If $(e_{1},\dots,e_{n})$ is
a local basis of vector fields on $U$
and $(\epsilon^{1},\dots,\epsilon^{n})$ is the
dual basis, then the $1$-form-valued vector
$\nabla X$ in $U$
can be written as $\sum_{j=1}^n \nabla_{e_{j}}X \otimes
\epsilon^{j} \ .$ Thus $i_{\nabla X} = \nabla_{e_{j}}X  \wedge
i_{\epsilon^{j}} \ .$
Therefore, by Proposition \ref{propdivconn} and Lemma \ref{OOO},
$$\begin{array}{rcl}
\ddiv_\nnabla(i_{\nabla X})
&=& \ddiv_\nnabla(\sum_{j=1}^n
\nabla_{e_{j}}X\wedge i_{\epsilon^{j}})\\
 &=&\sum_{j=1}^n
\nabla_{e_{j}}X\ddiv_\nnabla(i_{\epsilon^{j}}) -
\sum_{j=1}^n
\epsilon^{j}(\nabla_{e_{j}}X)= - \ddiv_\nabla(X) \ .
\end{array}
$$
It follows that
$$
\Delta^{Schouten,\nnabla}(X)
= - {\frac 12}\ddiv_\nnabla(\lcf X, ~\cdot~\rcf)
= - \ddiv_{\nabla}(X) \ .
$$
It follows from Corollaries \ref{dercurv}
and \ref{vvv} together with Proposition \ref{curvtr} that,
if $\nabla$ is flat, the operator $(\Delta^{Schouten,\nnabla})^{2}$ is a
derivation of the sheaf of multivectors with respect to the exterior
product. Since, moreover,
$(\Delta^{Schouten, \nnabla})^{2}$ is of $\mathbb Z$-degree $- 2$, it
vanishes.
\qed

\medskip

\noindent{\it Remark.}
Koszul \cite{Kz1} proves
that, conversely, any generator of the Schouten bracket 
of multivectors on a manifold $M$
is of the form $\Delta ^{\nabla}$ for some 
torsionless linear connection on $M$, and that two connections 
give rise to the same generator of the
Schouten bracket if and only if they  induce
the same linear
connection on $\bigwedge^m TM$, $m$ being the dimension of the manifold.  
We have not found any straightforward extension of this result to the case of
odd Poisson brackets on supermanifolds in general. 

\subsubsection{Conclusion}

Given a smooth manifold $M$, we set $A=C^{\infty}(M)$, and we let
${\rm Der}~A$ denote the module of vector fields on $M$.
The definition of a divergence operator on a graded algebra 
reduces, in the purely even case of $A$, to the requirement that
the linear operator, 
$\ddiv : {\rm Der}~A \to A$, satisfy the identity $\ddiv(fX) = f\ddiv X
+ X(f)$, for any $f \in A$ and $X \in {\rm Der}~A$.
The operators $\ddiv_{\nabla}$, considered in Section \ref{gensch},
where $\nabla$ are torsionless
linear connections on $M$, 
are examples of divergence operators. 
Other examples are furnished by the operators
$\ddiv_{\mu}$
associated to volume forms, ${\mu}$, on an orientable manifold $M$.
The Schouten bracket is an odd Poisson bracket on
the graded commutative, associative algebra, 
${\bf \Lambda} = \bigwedge _A ({\rm Der}~A)
= \oplus_{k=0}^{m} \bigwedge ^k _A ({\rm Der}~A)$, where
$m$ is the dimension of the manifold $M$.
It is the opposite of a divergence operator that can be extended into a
generator of the Schouten bracket.
In fact, for any divergence operator on $A$, 
the operator $- \ddiv$ can be uniquely extended to a generator
of ${\mathbb Z}$-degree $-1$, 
denoted 
$\Delta ^{(\ddiv)}$, of the Schouten bracket.
One can characterize the generator $\Delta ^{(\ddiv)}$
recursively since, for any $f \in A$, 
it commutes with the interior product $i_{\d f}$. 
More generally, for any form $\alpha$, 
$$
[i_{\alpha} , 
\Delta ^{(\ddiv)}]= - i~_{\d\alpha} \ .
$$

It is easy to see that, in the purely even case, a divergence
operator $\ddiv : \Der A \to A$ is 
nothing but a right $(A,\Der A)$-connection on $A$, in the sense of
Huebschmann \cite{Hu} \cite{Hu2}. In fact,
if $\ddiv$ is a divergence operator, then 
$$
(f,X) \in A \times \Der A
\mapsto f \circ X = - \ddiv(fX) \in A
$$ 
is a right $(A,\Der A)$-connection on
$A$, and, conversely, if  $(f,X) \in A \times \Der A
\mapsto f \circ X \in A$ is a right $(A,\Der A)$-connection on $A$, then
$X \in \Der A \mapsto - 1 \circ X$, where $1$ is the unit of $A$, is a
divergence operator. Moreover,
the right $(A,\Der A)$-connection is a right $(A,\Der A)$-module
structure if and only if the curvature of the divergence operator, defined
by \eqref{curvdiv}, vanishes. In fact,
$$
(f \circ X_1)\circ X_2 - (f \circ X_2) \circ X_1 - f \circ [X_1,X_2] =
f {\mathcal R}^{\ddiv}(X_1,X_2) \ .
$$
In the papers cited above, 
the notion of a divergence operator does
not appear explictly, but the preceding remarks show that 
the 1-to-1 correspondence (\cite{Hu}, Theorem 1) between 
right $(A,\Der A)$-connections on $A$ and generating operators
of the Schouten bracket of ${\bf \Lambda}$ 
yields a
1-to-1 correspondence between divergence operators and generators, which
restricts to a 1-to-1 correspondence between divergence operators whose
curvature vanishes and generators whose square vanishes.
Also, the 1-to-1 correspondence (\cite{Hu}, Theorem 3)
between right $(A,\Der A)$-connections on
$A$ and left $(A,\Der A)$-connections on the top exterior power,
$\bigwedge ^m _A \Der A =\bigwedge^{\rm top}_A \Der A $,
translates into a 1-to-1 correspondence between divergence operators and
left $(A,\Der A)$-connections on the top exterior power. 
The canonical bundle 
$\bigwedge^{\rm top}_A \Der A$ is, in a natural way, a right module;
equipping it with a left module structure, which can be done by
choosing a volume element, is equivalent to equipping 
$A$ itself with a
right module structure and therefore to selecting a divergence
operator
whose curvature vanishes (cf Propositon \ref{curvber}).
To summarize, divergence operators,
right connections on $A$, left connections
on the top exterior power of $\Der A$, 
and generators of the Schouten bracket are in 1-to-1 correspondence.
The definition of divergence operators and the preceding
constructions extend to the framework of Lie algebroids and to that of
Lie-Rinehart algebras.
See \cite{Hu}, \cite{Hu2} and also \cite{X} and \cite{KS}.
While there is a functor from 
Lie-Rinehart algebras to Gerstenhaber algebras, there is also a
functor from
Lie-Rinehart algebras with a divergence operator (resp., divergence
operator with
vanishing curvature) to Gerstenhaber algebras with a generator (resp.,
to Batalin-Vilkovisky algebras).
In the case of a complex analytic manifold $M$ and its algebra of
analytic functions, 
the left $(A,\Der A)$-module structures on the canonical 
bundle (top
exterior power of holomorphic vector fields) 
are called Calabi-Yau structures \cite{Sch}. In this case, 
left (resp., right) $(A, \Der A)$-module structures coincide with
left (resp., right) ${\mathcal D}_M$-module structures.

The extension 
of the above 1-to-1 correspondences valid in the purely even case 
to the case where $A$ itself is a
$\mathbb Z$- or ${\mathbb Z}_{2}$-graded algebra, ${\bf A}$, 
remains to be done. 
The appropriate framework is that of the graded Lie-Rinehart algebras, whose
theory has already been developped by Huebschmann (1990, unpublished),
and left and right $({\bf A}, \Der{\bf A})$-connections and module
structures
in an appropriate sense. Sheaves of graded Lie-Rinehart algebras
should then be considered, the fundamental example being 
$({\mathcal A}, \Der{\mathcal A})$, for
any supermanifold $(M,{\mathcal A})$.
A divergence operator with vanishing curvature should define a 
right $({\mathcal A}, \Der{\mathcal A})$-module structure on
${\mathcal A}$, 
and there should be 1-to-1 correspondences between divergence
operators,
right structures on the structural sheaf and left
structures on the berezinian sheaf. 
Another approach is by means of the theory of
${\mathcal D}$-modules. Left and right
${\mathcal D}$-modules on complex
supermanifolds have been studied by Penkov \cite{P}, who showed that the
berezinian sheaf of a complex analytic supermanifold is a right 
${\mathcal D}$-module in
a canonical way.
Defining a left ${\mathcal D}$-module structure on the berezinian
sheaf, which can be done by choosing a berezinian volume, is
equivalent to defining a right ${\mathcal D}$-module structure on the
structural sheaf, and should be equivalent to
the choice of a divergence operator.

One can define graded analogues of the modules of multivectors on a
manifold as modules of skew-symmetric multiderivations
of ${\bf A}$, and one can generalize these notions to the case of
sheaves of graded algebras over a manifold.
Multigraded generalizations of the
Schouten bracket on the space of skew-symmetric
multiderivations of a graded algebra
were defined by Krasil'shchik in
\cite{Kr2}, following his earlier paper \cite{Kr1}.
An analogue of the 1-to-1 correspondence between
divergence operators and generators 
should be also valid in the graded case.

\noindent{\it Conjecture. 
A divergence operator on the graded algebra ${\bf A}$,
up to sign factors, can be uniquely extended to 
an operator on the skew-symmetric
multiderivations of ${\bf A}$
that generates, in a suitable sense, the bigraded 
Krasil'shchik-Schouten bracket.} 

In particular, this construction would associate to a 
divergence operator
on a supermanifold $(M,\mathcal A)$ a 
generator of the bigraded bracket on multivectors on the
supermanifold.
We hope to return to this question and also to study the
relationship between the generators of a graded
bracket and those of its derived brackets, in the sense of \cite{KSf},
in a future publication.

\section*{Appendix. The berezinian sheaf}

We shall recall the definition of the berezinian integral and some
fundamental results, following
\cite{Be}, \cite{Ma}, \cite{Ro}, \cite{Vor} 
and mostly \cite{HM1} and \cite{HM2}.

Let $(M,{\mathcal A})$ be a supermanifold of dimension $m|n$, in the sense
of \cite{Ko}. Thus, $M$ is a smooth manifold and ${\mathcal A}$ is a 
sheaf of
${{\mathbb Z}}_2$-graded commutative, associative ${\mathbb R}$-algebras
over $M$. There is an exact sequence
$$0\to {\mathcal N}\to {\mathcal A}\to {\mathcal A}/{\mathcal N}\to 0 \ 
,$$
where ${\mathcal N}$ is the sheaf of nilpotent sections of ${\mathcal A}$,
and ${\mathcal
A}/{\mathcal N}$ is the sheaf $C^\infty(M)$, regarded as trivially
graded.
The projection ${\mathcal A}\to {\mathcal A}/{\mathcal N}=C^\infty(M)$ is 
denoted
by the symbol \ $\widetilde{~~~} \ ,$
and there is a unique prolongation
to
the module of differential forms of this projection,
that commutes with the de~Rham differentials. Thus, if we denote by $\d$
and
$\d^G$ the de~Rham differentials in $M$ and $(M,{\mathcal A})$, then
$$\widetilde{\d^G\alpha} = \d\widetilde\alpha,$$
for any differential form $\alpha$ on the supermanifold $(M,{\mathcal 
A})$.

The {\it berezinian sheaf} can be described as follows.
Let ${\mathcal P}^k({\mathcal A})$ be the vector space of the differential
operators
of order $k$ on ${\mathcal A}$. There
is both a right and a left ${\mathcal A}$-module structure on
${\mathcal P}^k({\mathcal A})$ given by
$(f.P)(g) = f.P(g)$ and
$(P.f)(g) = P(f.g)$, respectively, for sections  $f,g$ of ${\mathcal A}$ 
and
$P\in {\mathcal P}^k({\mathcal A})$.
If $(x^1,\dots,x^m,s^1,\dots,s^n)$
are graded coordinates on an open set $U$ in $M$,
then
${\mathcal P}^k({\mathcal A})|_{U}$ is free for both structures of 
${\mathcal
A}$-module,
with basis
$$(\frac {\partial}{\partial x^{1}})^{k_1} \circ\dots\circ
(\frac {\partial}{\partial x^{m}})^{k_{m}}\circ
\frac {\partial}{\partial s^{\rho_1}}\circ\dots\circ
\frac {\partial}{\partial s^{\rho_j}},
$$
where $k_1, \ldots, k_{m} \in{\mathbb N}$, $1 \leq \rho_{1} < \rho_{2} < 
\ldots
< \rho_{j} \leq n $ and
$k_1+\dots +k_m+j = k$.

Let us now define ${\mathcal P}^k({\mathcal A},\Omega^m_{\mathcal A}) =
\Omega^m_{\mathcal A}\otimes {\mathcal P}^k({\mathcal A})$, where 
$\Omega^m_{\mathcal
A}$
is the sheaf of differential $m$-forms on $(M,{\mathcal A})$.
Let ${\mathcal K}^{n}$
be the subsheaf of elements $\bf P$ in ${\mathcal P}^n({\mathcal 
A},\Omega^m_{\mathcal A})$
such that, for any
section $f$ of ${\mathcal A}$ over an open set $U$ of $M$ with
compact support, there exists an $(m-1)$-differential form
$\omega$ with compact support in $U$
such that $\widetilde{{\bf P}(f)} = \d\omega$.
Is is easy to show that ${\mathcal K}^{n}$ is a subsheaf of right 
${\mathcal
A}$-modules of
${\mathcal P}^n({\mathcal A},\Omega^m_{\mathcal A})$.
The berezinian sheaf is the quotient sheaf
${\mathcal P}^n({\mathcal A},\Omega^m_{\mathcal A}) /{\mathcal K}^{n}.$
The sections of this sheaf can be locally expressed as
$$[ \d^Gx^{1}\wedge\dots\wedge\d^Gx^{m} \otimes
{\frac {\partial}{\partial s^{1}}} \circ \dots \circ
{\frac {\partial}{\partial s^{n}}} ].f \ ,$$
where $f$ is a section of ${\mathcal A}$.
If $V\subset M$ is an open set with
graded coordinates
$(y^1,\dots,y^m,t^1,\dots,t^n)$, then, on $U \cap V$,
$$
[\d^Gy^1\wedge\dots\wedge\d^Gy^m \otimes
\frac\partial{\partial t^1}\circ\dots\circ\frac\partial {\partial {t^n}}]
=
[\d^Gx^1\wedge\dots\wedge\d^Gx^m \otimes
\frac\partial{\partial s^1}\circ\dots\circ\frac\partial
{\partial {s^n}}].Ber\begin{pmatrix}  A&B\\ C&D\end{pmatrix}$$
where
$$\begin{pmatrix}  A&B\\ C&D\end{pmatrix} =
\begin{pmatrix}
\big(\frac {\partial y^{i}}{\partial x^{j}}\big)&
\big(\frac {\partial t^{\rho}}{\partial x^{j}}\big)\\
\big(\frac {\partial y^{i}}{\partial s^{\sigma}}\big)&
\big(\frac {\partial t^{\rho}}{\partial s^{\sigma}}\big)
\end{pmatrix}
$$
and where $Ber$ denotes the berezinian.
(The berezinian, or superdeterminant, of an invertible even matrix of the 
form $\begin{pmatrix} A&B\\
C&D \end{pmatrix}$
is $\det(A-BD^{-1}C)\det(D^{-1})$.)

If $M$ is an orientable smooth manifold,
the {\it Berezin integral},~ $\int_{(M,{\mathcal A})}$ \ , maps the 
sections
with compact
support of the berezinian sheaf to $\mathbb R$, and is defined by
$\int_{(M,{\mathcal A})} [{\bf P}] = \int_M \widetilde{{\bf P}(1)}.$
As an example, if $(M,{\mathcal A}) = {\mathbb R}^{m|n}$, then
$$\int_{{\mathbb R}^{m|n}}
[\d^Gx^1\wedge\dots\wedge\d^Gx^m \otimes
\frac\partial{\partial s^1}\circ\dots\circ\frac\partial {\partial {s^n}}]
.f = \hskip -3pt
(-1)^{{\frac {n(n-1)}{2}}} \int_{{\mathbb R}^n} f_{(1,2,\dots,n)}(x) \d
x^{1}\wedge
\dots\wedge\d x^{m},$$
where $f_{(1,2,\ldots,n)}$ is the coefficient of $s^{1}s^{2}\ldots s^{n}$
in the expansion of $f$ as a sum of products of the $s^{\rho}$'s.

A section, $\xi$, of the berezinian sheaf is called a {\it berezinian
volume} if
it is a generator of the berezinian sheaf, {\it i.e.}, if any other 
section can
be uniquely written as $\xi.f$ for some section $f$ of ${\mathcal A}$.
A berezinian volume is a homogeneous section of the
berezinian sheaf, whose degree depends on the parity of the dimension $n$.
If $\xi$ is a berezinian volume and $v$ is a section of ${\mathcal A}$,
then $\xi.v$ is
also a berezinian volume if and only $v$ is invertible and even.

In order to define the {\it Lie derivatives} of berezinian volumes with 
respect
to graded vector fields, we first observe that,
in a similar way, we can define the right submodule ${\mathcal K}^{n+k}$ 
of
${\mathcal P}^{n+k}({\mathcal A}, \Omega_{\mathcal A}^m)$, for each $k 
\ge 1$, and
that the canonical inclusion
${\mathcal P}^{n}({\mathcal A}, \Omega_{\mathcal A}^m)\hookrightarrow
{\mathcal P}^{n+k}({\mathcal A}, \Omega_{\mathcal A}^m)$ induces an 
isomorphism of
sheaves of right ${\mathcal A}$-modules from
${\mathcal P}^{n}({\mathcal A}, \Omega_{\mathcal A}^m)/{{\mathcal 
K}^{n}}$ to
${\mathcal P}^{n+k}({\mathcal A}, \Omega_{\mathcal A}^m)/{\mathcal 
K}^{n+k}.$

Let $D$ be a graded vector fied on $(M,{\mathcal A})$.
The Lie derivative of the berezinian volume $[\omega \otimes P]$
with respect to $D$ is
$${\mathcal L}_D[\omega\otimes P] = (-1)^{|D||\omega\otimes P|+1}
[\omega\otimes P\circ D].
$$
The main properties of the Lie derivatives of berezinian volumes are
stated in Section \ref{2}, 
and are used there in order to derive the properties
of the divergence operators.

\bigskip

\noindent{\bf Acknowledgments}. {\small 
Y. K.-S. would like to thank
P. Cartier, J. Huebschmann, D. Leites, H. Khudaverdian, A. Schwarz,
J. Stasheff and T. Voronov for interesting exchanges
on the topic of this work. Thanks are also due to the referee who
identified various weaknesses in the first version of this paper.

J. M. is partially
supported by Pla Valenci{\`a} de
Ci{\`e}ncia i Tecnologia, grant $\#POST99$-$01$-$30$, and
by DGICYT, grant $\#PB97$-$1386$.}


\begin{thebibliography}{9999}


\bibitem{ASZK} Alexandrov, M.,
Kontsevich, M., Schwarz, A., Zaboronsky, O.,
The geometry of the master equation and topological quantum field
theory, {\sl Int. J. Mod. Phys.} {\bf A12} (1997), 1405-1429.


\bibitem{B} Batalin, I. A., Vilkovisky, G. A., Gauge algebra and
quantization,
{\sl Phys. Lett.}
{\bf B 102} (1981), 27-31.

\bibitem{BV} Batalin, I. A., Vilkovisky, G. A., Closure of the 
gauge algebra, generalized Lie equations and Feynman rules, Nuclear
Physics {\bf B 234} (1984), 106-124.

\bibitem{Be} Berezin, F. A., {\sl Introduction to Superanalysis},
D. Reidel (1987).

\bibitem{BM} Beltr{\'a}n, J. V., Monterde, J.,
Graded Poisson structures on the algebra of differential forms,
{\sl Comment. Math. Helv.}  {\bf 70} (1995), 383-402.

\bibitem{BMS-V} Beltr{\'a}n, J. V., Monterde, J., 
S{\'a}nchez-Valenzuela, O. A.,
Graded Jacobi operators on the algebra of differential forms,
{\sl Compositio Math.} {\bf 106} (1997), 43-59.

\bibitem{D}
P. Deligne {\sl et al.}, eds., {\sl Quantum Fields and Strings: A Course
for
Mathematicians}, Amer. Math. Soc. (1999), vol. 1, part 1.

\bibitem{dW} DeWitt, B., {\sl Supermanifolds}, Cambridge Univ. Press (1984).

\bibitem{FN}
Fr{\"o}licher, A., Nijenhuis, A., Theory of vector-valued differential
forms, part I, {\sl Indag. Math.}, {\bf 18} (1956), 338-359.

\bibitem{G} Getzler, E., Batalin-Vilkovisky algebras and two-dimensional
topological field theories, {\sl Commun. Math. Phys.} {\bf 159} (1994),
265-285.

\bibitem{Z} Hata, H., Zwiebach, B.,
Developing the covariant Batalin-Vilkovisky approach to string theory,
{\sl Ann. Phys.} {\bf 229} (1994), 177-216.

\bibitem{HM1} Hern{\'a}ndez Ruip{\'e}rez, D., Mu{\~n}oz Masqu{\'e}, J.,
Construction intrins{\`e}que du faisceau de Berezin d'une vari{\'e}t{\'e}
gradu{\'e}e, {\sl Comptes Rendus Acad. Sci. Paris, S{\'e}r. I Math.}
{\bf 301} (1985), 915-918.

\bibitem{HM2} Hern{\'a}ndez Ruip{\'e}rez, D., Mu{\~n}oz Masqu{\'e}, J.,
Variational Berezinian problems and their relationship with graded
variational problems, {\sl Diff. Geometric Methods in Math. Phys.
(Salamanca 1985),
Lect. Notes  Math.} {\bf 1251}, Springer-Verlag (1987), 137-149.


\bibitem{Hu1} Huebschmann, J., Poisson cohomology and quantization,
{\sl J. f{\"u}r die reine und angew. Math.} {\bf 408} (1990), 57-113.

\bibitem{Hu} Huebschmann, J., Lie-Rinehart algebras, Gerstenhaber
algebras, and Batalin-Vilkovisky algebras,
{\sl Ann. Inst. Fourier} {\bf 48} (1998), 425-440.

\bibitem{Hu2} Huebschmann, J., Duality for Lie-Rinehart algebras
and the modular class, 
{\sl J. f{\"u}r die reine und angew. Math.} {\bf 510} (1999), 103-159.

\bibitem{kh0}
Khudaverdian, O. M., Geometry of superspace with even and odd
brackets, {\sl J. Math. Phys.} {\bf 32} (1991), 1934-1937.


\bibitem{kh}
Khudaverdian, O. M.,
Batalin-Vilkovisky formalism and odd symplectic geometry, {\sl
Geometry and integrable models (Dubna 1994)}, Pyatov, P. N., 
Solodukhin, S. N., eds.,
World Sci. Publish. (1996), 144--181.



\bibitem{khn}
Khudaverdian, O. M., Nersessian, A. P., On the geometry of the
Batalin-Vilkovisky formalism, {\sl Mod. Phys. Lett.} {\bf A 8} (1993),
2377-2385.


\bibitem{KSf}
Kosmann-Schwarzbach, Y., From Poisson algebras 
to Gerstenhaber algebras, {\sl Ann. Inst. Fourier}
{\bf 46} (1996), 1243-1274.

\bibitem{KS}
Kosmann-Schwarzbach, Y., Modular vector fields and Batalin-Vilkovisky
algebras,
{\sl Banach Center Publications} {\bf 51} (2000), 109-129.

\bibitem{KSM}
Kosmann-Schwarzbach, Y., Magri, F., Poisson-Nijenhuis structures,
{\sl Ann. Inst. Henri Poincar{\'e}} {\bf A53} (1990), 35-81.

\bibitem{Ko}
Kostant, B., Graded manifolds, graded Lie theory and prequantization,
{\sl Proc. Conf. Diff. Geom. Methods in Math. Phys. (Bonn
1975), Lecture Notes  Math.}
{\bf 570}, Springer-Verlag (1977), 177-306.


\bibitem{Kz1} Koszul, J.-L.,
Crochet de Schouten-Nijenhuis et cohomologie,
in {\sl {\'E}lie Cartan et les
math{\'e}\-ma\-tiques d'aujourd'hui},
Ast{\'e}risque, hors s{\'e}rie, Soc. Math. Fr.
(1985), 257-271.

\bibitem{Kr1} Krasil'shchik, I. S., Schouten brackets and canonical
algebras, {\sl Lecture Notes Math.} {\bf 1334}, Springer-Verlag
(1988), 79-110.

\bibitem{Kr2}  Krasil'shchik, I. S., Supercanonical algebras and
Schouten
brackets, {\sl Mat. Zametki}, {\bf 49}(1) (1991), 70-76, {\sl
Mathematical Notes} {\bf 49}(1) (1991), 50-54.

\bibitem{L0} Leites, D.,
{\sl Supermanifold Theory}, 
Karelia Branch of the USSR Acad. of Sci., Petrozavodsk,
(1983) (in Russian).

\bibitem{L} Leites, D., Quantization and supermanifolds, Supplement 3 in
Berezin, F. A., Shubin, M. A., {\sl The Schr{\"o}dinger Equation}, Kluwer
(1991).


\bibitem{LZ} Lian, B. H., Zuckerman, G. J., New perspectives on the
BRST-algebraic structure of string theory, {\sl Commun. Math. Phys.}
{\bf  154} (1993), 613-646. 

\bibitem{Ma} Manin, Y. I., {\sl Gauge Field Theory and Complex Geometry},
Springer-Verlag (1988).


\bibitem{MP} 
Manin, Y. I., Penkov, I. B., The formalism of left and 
right connections
on supermanifolds,
{\sl Lectures on Supermanifolds, Geometrical Methods and Conformal Groups},
Doebner, H.-D., Hennig, J. D., Palev, T. D., eds., World
Sci. Publ. (1989), 3-13.

\bibitem{MM} Monterde, J., Montesinos, A., Integral curves of derivations,
{\sl Ann. Global Anal. Geom.} {\bf 6} (1988), 177-189.

\bibitem{MS} Monterde, J., S{\'a}nchez Valenzuela, A.,
The exterior derivative as a Killing vector field,
{\sl Israel J. Math.} {\bf 93} (1996), 157--170.

\bibitem{P} Penkov, I. B., ${\mathcal D}$-modules 
on supermanifolds, {\sl Invent. Math.} {\bf 71} (1983), 501-512.

\bibitem{Ro} Rothstein, M., Integration on noncompact supermanifolds,
{\sl Trans. Amer. Math. Soc.} {\bf 299} (1987), 387-396.

\bibitem{Sch} Schechtman, V., {\sl Remarks on formal deformations and
Batalin-Vilkovisky algebras}, preprint math.AG/9802006.
 
\bibitem{Sc} Schwarz, A., Geometry of Batalin-Vilkovisky quantization,
{\sl Commun. Math. Phys.} {\bf 155} (1993), 249--260.

\bibitem{Sc2} Schwarz, A., Semi-classical approximation in
Batalin-Vilkovisky
formalism,
{\sl Commun. Math. Phys.} {\bf 158} (1993), 373-396.

\bibitem{St} Stasheff, J., Deformation theory and the 
Batalin-Vilkovisky master equation, {\sl Deformation Theory
and Symplectic Geometry (Ascona 1996)}, Sternheimer, D., Rawnsley, J.,
Gutt, S., eds., Kluwer (1997), 271-284.

\bibitem{V94} Vaisman, I., {\sl Lectures on the Geometry of Poisson
Manifolds}, Birkh{\"a}user, 1994.

\bibitem{Vor} Voronov, T., Geometric integration theory
on supermanifolds,
{\sl Sov. Sci. Rev.}, C Math., {\bf 9} (1992), 1-138.

\bibitem{We} Weinstein, A.,
The modular automorphism group of a Poisson manifold,
{\sl J. Geom. Phys.} {\bf 23} (1997), 379-394.

\bibitem{W} Witten, E., A note on the antibracket formalism, {\sl
Mod. Phys. Lett.} {\bf A5} (1990), 487-494.

\bibitem{X} Xu, P., Gerstenhaber algebras and BV-algebras in Poisson
geometry, {\sl Commun. Math. Phys.} {\bf 200} (1999), 545-560.

\end{thebibliography}
\end{document}